\documentclass{amsart}
\usepackage{amssymb,stmaryrd}
\usepackage{amsfonts}
\usepackage{amstext}
\usepackage{algorithmic}
\usepackage{algorithm}
\usepackage{graphicx}
\usepackage{epstopdf}
\usepackage[all]{xy}
\usepackage{MnSymbol,slashed}

\numberwithin{equation}{section}
\newtheorem{theorem}{Theorem}[section]
\newtheorem{lemma}[theorem]{Lemma}
\newtheorem{proposition}[theorem]{Proposition}
\newtheorem{corollary}[theorem]{Corollary}

\theoremstyle{definition}

\newtheorem{example}[theorem]{Example}

\theoremstyle{remark}
\newtheorem{remark}[theorem]{\bf{Remark}}

\renewcommand{\S}{{\mathcal{S}}}
\newcommand{\R}{{\mathbb{R}}}

\newcommand{\C}{{\mathbb{C}}}
\newcommand{\Z}{{\mathbb{Z}}}
\newcommand{\T}{{\mathbb{T}}}

\newcommand{\<}{{\langle}}

\renewcommand{\)}{{)}}
\renewcommand{\>}{{\rangle}}

\newcommand{\CS}{{\mathcal{S}}}

\newcommand{\CJ}{{\mathcal{J}}}
\newcommand{\wedgeq}{{\wedge\kern-5pt\cdot}}

\newcommand{\tens}{\otimes}

\newcommand{\id}{{\rm id}}

\newcommand{\extd}{{\rm d}}
\newcommand{\del}{{\partial}}
\newcommand{\eps}{\epsilon}

\newcommand{\la}{{\triangleright}}

\newcommand{\dirac}{{ \slashed{D} }}

\newcommand{\doublenabla}{%
  \nabla\mkern-12mu\nabla}   
  
\begin{document}

\title{Geometric Dirac operator on noncommutative torus and $M_2(\C)$}


\author{E. Lira-Torres and S. Majid}
\address{Queen Mary, University of London\\
School of Mathematics, Mile End Rd, London E1 4NS, UK}
\thanks{The first author was partially supported by CONACyT and Fundaci\'on
  Alberto y Dolores Andrade (M\'exico)}

\email{s.majid@qmul.ac.uk}


\begin{abstract}  We solve for quantum-geometrically realised spectral triples or `Dirac operators' on the noncommutative torus $\C_\theta[T^2]$ and on the algebra  $M_2(\C)$ of $2\times 2$ matrices with their standard quantum metrics and associated quantum Levi-Civita connections. For  $\C_\theta[T^2]$, we obtain a even standard spectral triple but now uniquely determined by full geometric realisability. For $M_2(\C)$, we are forced to the flat quantum Levi-Civita connection and again obtain a natural fully geometrically realised even spectral triple. In both case there is also an odd spectral triple for a different choice of a sign parameter. We also consider an alternate quantum metric on $M_2(\C)$ with curved quantum Levi-Civita connection and find a natural 2-parameter of almost spectral triple in that $\dirac$ fails to be antihermitian. In all cases, we split the construction into a local tensorial level related to the quantum geometry, where we classify the results more broadly, and the further requirements relating to the Hilbert space structure. We also illustrate the Lichnerowicz formula for $\dirac^2$ which applies in the case of a  full  geometric realisation.  \end{abstract}
\maketitle 

\section{Introduction}

Noncommutative geometry, the idea that coordinates could be noncommutative as in quantum mechanics, has gained interest in recent years as a plausibly better description of spacetime that accounts for quantum gravity effects\cite{Sny,Ma:pla,MaRue,DFR,Hoo}, as well as clearly having potential for diverse applications in quantum theory. Probably the best known approach is the one coming out of cyclic cohomology and operator algebras\cite{Con,Bon} with the noncommutative torus as a prime example\cite{Con:dif} and leading up to the notion of `Connes spectral triple' in the role of Dirac operator\cite{Con0,ConMar}. Meanwhile, since the arrival of true quantum groups in the 1980s, there has emerged a constructive `quantum groups' approach which starts with a unital but possibly noncommutative `coordinate algebra' $A$ with quantum or classical symmetries to narrow down the choice of a bimodule of 1-forms $\Omega^1$. In recent years, this was extended to a systematic theory of  `quantum Riemannian geometry' (QRG)\cite{BegMa} with a quantum metric $g\in \Omega^1\tens_A\Omega^1$, a quantum Levi-Civita connection (QLC) $\nabla:\Omega^1\to \Omega^1\tens_A\Omega^1$, etc. This framework of curved quantum spacetimes has been applied to baby models of quantum gravity and quantum geodesics, and is also beginning to be applied to particle physics, see \cite{ArgMa4} for a review. 

Of interest in this paper and linking the two approaches to noncommutative geometry is a framework \cite{BegMa:spe}\cite[Chapter~8.5]{BegMa} to extend a QRG  with (a) a `spinor' bimodule $\CS$ equipped with a bimodule connection $\nabla_\CS$, (b) a `Clifford action' bimodule map $\la:\Omega^1\tens_A\CS\to \CS$ and (c) a `charge conjugation' bimodule map $\CJ:\CS\to \bar \CS$ to the conjugate bimodule (so $\CJ$ is antilinear), so as to arrive at a quantum-geometric  Dirac operator 
\[ \dirac=\la\circ\nabla_\CS:\CS\to \CS\]
 obeying the `local tensorial' part of Connes axioms for a spectral triple (by which we mean before consideration of a Hilbert space structure on $\CS$). If the latter step can also be completed then we obtain a `quantum-geometrically realised' spectral triple. Not every Connes spectral triple can be geometrically realised like this but, rather, geometric realisability, picks out natural spectral triples where we have this additional underlying quantum geometry.  Moreover, \cite{BegMa:spe} proposed two further (optional) restrictions on the Clifford action $\la$ which are not needed for a spectral triple but  link the spinor structures to the metric as would be the case classically. The most important is that
$\la$ should be covariantly constant, i.e. covariant with respect to the connections on the two sides of the map, see (\ref{covla}). The second is  `Clifford relations', see (\ref{cliffla}), which we allow up to an automorphism as in \cite{BegMa:spe}. If both of these hold, we will say that $\dirac$ has a `full' geometrically realisation. 

In previous work, \cite{BegMa:spe} constructed a natural $q$-deformed Dirac operator on the standard $q$-sphere $\C_q[S^2]$ as a full geometric realisation at the local tensorial level except with $\CJ$ being not quite an antilinear isometry (but a $q$-deformed one). Similarly, \cite{LirMa2} analysed the fuzzy sphere and found a natural fully geometrically realised spectral triple as  a member of a previous family proposed in \cite{And}.  The work \cite{BegMa:spe}  also gave an example of a quantum-geometrically realised spectral triple on $M_2(\C)$ but  without (\ref{covla}) for covariance of $\la$, so not a full one. These results, as well as our new results described below, all have 2-dimensional spinor bundles and are complementary to a recent canonical construction from a QRG\cite{Ma:dir}, which would be 3-dimensional in these examples. That general construction is of interest as a version of the Hodge-Dirac operator rather than of a conventional spinor Dirac operator and does not solve (\ref{cliffla}), so again not a full geometric realisation.

We now continue this geometric realisation programme with, in Section~\ref{sector}, a similar analysis for the algebraic noncommutative torus $\C_\theta[\T^2]$ where we  obtain the `obvious' construction
\[ (\dirac \psi)_\beta= \del_i\psi_\alpha \sigma^{i\alpha}{}_\beta.\]
Here $\sigma^i$ are Pauli matrices for $i=1,2$, $\del_i$ are certain natural derivations\cite{Con:dif} and $\psi_\alpha$ is a 2-spinor with values in the noncommutative torus. Up to conventions, this is also the Euclidean metric case in the exposition~\cite{CPR}. Our result is that this $\dirac$ is unique up to unitary equivalence among fully quantum-geometrically realised spectral triples with $\CS$ free and 2-dimensional over the algebra. In the analysis, we allow  $\nabla$ be any weak QLC, but the construction forces this to be the unique QLC on this algebra given simply by $\nabla s^i=0$ on a canonical left-invariant self-adjoint basis $\{s^i\}$, with respect to which we take the quantum metric $g=s^1\tens s^1+s^2\tens s^2$.  More precisely, we actually have two spectral triples which one could view as $n=1,2$ respectively in the period 8 classification, ignoring the `even structure' $\gamma$ in the odd case if we want. It is known that there are four (equivariant) spectral triples on the noncommutative torus\cite{PasSit1} and presumably we obtain two of them. We then show in Section~\ref{sectorext}  that we get a much bigger moduli of $\nabla_\CS,\nabla$ if we drop (\ref{covla}) for covariance of $\la$ and just seek a geometric realisation, not necessarily full. Remarkably, the final $\dirac$ itself remains of the above form but with the possible addition of a constant mass term or chiral `mass' term.   

 Sections~\ref{secM2}  then provides a similar analysis for $M_2(\C)$, using the standard quantum geometry in \cite[Example~8.13]{BegMa}. We convert this to a self-adjoint basis $\{s^i\}$ for $\Omega^1$ where the standard quantum metric becomes
  \[ g=\imath (s^2\tens s^1-s^1\tens s^2),\]
  and we work with a natural 1-parameter sub-moduli of the known 4-parameter QLCs for this. We then find, Theorem~\ref{thm_alt}, a generic  construction for a large class of spectral triples at the local tensorial level of $\nabla_\CS,\la,\CJ$ provided we are at the unique flat QLC in this class. With some simplifying assumptions, we are led to a solution (Proposition~\ref{propspecM2})  for which $\dirac$ is furthermore compatible with the usual inner product given by ${\rm Tr}:M_2(\C)\to \C$  on the spinor components to give $n=0,1$  fully geometrically realised spectral triples (depending on a sign $\eps'=\pm1$ and ignoring the `even structure' $\gamma$ in the odd case if we want). The `Clifford relations'  (with a nontrivial automorphism) hold and $\la$ is covariantly constant, making this much more canonical than a  previous geometrically realised spectral triple on $M_2(\C)$ in \cite{BegMa:spe}, where covariance of $\la$ did not hold.  After a unitary change of basis in Corollary~\ref{corM2dirac}, we obtain 
 \[ (\dirac \psi)_\beta={\imath\over 2 \sqrt{2}}\left( (\sigma^1\psi_\alpha+\psi_\alpha \sigma^2)\sigma^1{}^\alpha{}_\beta+\eps'(\sigma^2\psi_\alpha+\psi_\alpha\sigma^1)\sigma^3{}^\alpha{}_\beta \right)\]
as a natural fully geometrically realised spectral triple on $M_2(\C)$. 
  
For completeness, the appendix does the same analysis for another quantum metric \cite[Exercise 8.5]{BegMa}, which in our self-adjoint basis has the diagonal but `Lorentzian' form
  \[ g= - s^1\tens s^1+ s^2\tens s^2,\]
 where the - sign in the metric is needed for quantum symmetry. We find that this does not admit solutions for $\nabla_\CS,\la,\CJ$  if we impose the full Clifford relations (\ref{cliffrel}) but does admit solutions if we impose only half of these relations. These lead to a natural real 2-parameter family of proper geometrically realised almost spectral triples on $M_2(\C)$, see Theorem~\ref{M2evenspec} and Proposition~\ref{propM2dirac}. The two real parameters $s,t$ here fall on a circle $s^2+t^2=2$ and the QLC $\nabla$ is forced to have a specific value with curvature. These are `almost' spectral triples in that $\dirac$ turns out in the construction to be neither antihermitian nor hermitian, possibly attributable to the Lorentzian nature of the quantum metric. 
 
The  $M_2(\C)$ work involves  significant development of the previous analysis\cite{LirMa2}, done in Section~\ref{secspebasis}, so as to  cover central bases of $\Omega^1$ that need not obey the Grassmann algebra as differential forms, and to handle Christoffel symbols of $\nabla,\nabla_\CS$ that need not be $\C$-valued. This, and the required simplifications in the inner case in Section~\ref{secspeinn}, are some new general  results of the paper and reduce everything to solving nonlinear matrix equations, which we then solved using Mathematica.  Section~\ref{secele} provides  a very brief introduction to quantum Riemannian geometry as in \cite{BegMa} and to the geometric realisation programme  \cite{BegMa:spe}. Proposition~\ref{propLapS} and Lemma~\ref{lemRS} provide general formulae for the spinor Laplacian and spinor curvature operator under our assumptions of central bases. These enter into the  Lichnerowicz formula (\ref{diracsqmat}) for $\dirac^2$ which is known to apply\cite[Prop.~8.45]{BegMa} for any fully quantum-geometrically realised spectral triple, and Proposition~\ref{proplichM2} illustrates this for our new $\dirac$ on $M_2(\C)$ as a check on the calculations. The paper concludes in Section~\ref{secrem} with directions for further work. 
  
\section{Quantum geometric spinor formalism with central bases}\label{secpre} 

We start in  Section~\ref{secele} with a very  briefly outline the formalism of \cite{BegMa:spe}\cite[Chap~8.5]{BegMa} for the construction of Connes spectral triples from within `quantum Riemannian geometry' (QRG). Then, in Section~\ref{secspebasis}, we extend the case of a central basis $s^i$ of $\Omega^1$ in \cite{LirMa2} to now allow the $s^i$ to not be Grassmann and the Christoffel symbols of the quantum Levi-Civita connection to not be constants. We also analyse the Lichnerowicz formula for $\dirac^2$ in this setting. In Section~\ref{secspeinn}, we look at the special features when $\Omega^1$ is inner, which will be needed  for the matrix algebra $M_2(\C)$. 

\subsection{Elements of QRG}\label{secele} We work with $A$ a unital algebra, typically a $*$-algebra over $\C$, in the role of `coordinate algebra'. Differentials are formally introduced as a bimodule $\Omega^1$ of 1-forms equipped with a map $\extd:A\to \Omega^1$ obeying the Leibniz rule $\extd(ab)=(\extd a)b+a\extd b$. We assume this extends to an exterior algebra $(\Omega,\extd)$ with $\extd^2=0$ and $\extd$ obeying the graded-Leibniz rule, and where $A,\extd A$ generate $\Omega$ (this is more restrictive than a differential graded algebra). The exterior algebra is called {\em inner} if there exists $\theta\in \Omega^1$ such that $\extd =[\theta,\ \}$ is a graded derivation. We also use this term for $\Omega^1$ with $\extd=[\theta,\ ]$. 

A quantum metric is $g\in \Omega^1\tens_A\Omega^1$ and a bimodule map inverse $(\ ,\ ):\Omega^1\tens_A\Omega^1\to A$ in an obvious sense, together with some form of quantum symmetry condition. The general functorial definition of the latter is $\wedge(g)=0$, but in some contexts it  seems better to modify this. Following \cite{BegMa}, we say $g$ is a generalised quantum metric if no form of symmetry is imposed. A (left) bimodule connection\cite{DVMic,Mou} on $\Omega^1$ is $\nabla:\Omega^1\to \Omega^1\tens_A\Omega^1$ obeying 
\[ \nabla(a.\omega)=a.\nabla\omega+ \extd a\tens\omega,\quad \nabla(\omega.a)=(\nabla\omega).a+\sigma(\omega\tens\extd a)\]
for all $a\in A,\omega\in \Omega^1$, for some `generalised braiding' bimodule map $\sigma:\Omega^1\tens_A\Omega^1\to \Omega^1\tens_A\Omega^1$. The latter, if it exists, is uniquely determined and not additional data. But for an inner calculus, connections are given by a free choice of module map $\sigma$ and an additional bimodule map $\alpha:\Omega^1\to \Omega^1\tens_A\Omega^1$, see \cite[Prop.~8.11]{BegMa}. A connection is torsion free if $T_\nabla:=\wedge\nabla-\extd$ vanishes, and metric compatible if 
\[ \nabla g=(\nabla\tens \id+ (\sigma\tens\id)(\id\tens\nabla))g\]
vanishes. When both vanish, we have a {\em quantum Levi-Civita connection} (QLC). A weaker concept is a WQLC where we demand zero torsion and zero {\em cotorsion} in the sense
\[ (\extd\tens \id- \id\wedge \nabla)g=0\]
and we say in this case that we have a weak quantum Riemannian geometry (WQRG). The curvature of $\nabla$ is (similarly) defined as 
\[ R_\nabla=(\extd\tens \id- \id\wedge \nabla)\nabla:\Omega^1\to \Omega^2\tens_A\Omega^1\]
and the canonical Laplacian on $A$ is
\[ \square = (\ ,\ )\nabla\extd.\]
Finally, working over $\C$, we need  $(\Omega,\extd)$ to be a $*$-calculus, $g$ `real' in the sense $g^\dagger=g$, where $\dagger={\rm flip}(*\tens *)$ is well defined on $\Omega^1\tens_A\Omega^1$, and $\nabla$ $*$-preserving in the sense $\nabla\circ *=\sigma\circ\dagger\circ\nabla$, see\cite{BegMa}. 

The formulation of a bimodule connection above applies equally well to any $A$-bimodule  $\CS$ in the role of sections of a `vector bundle', now with $\nabla_\CS:\CS\to \Omega^1\tens_A\CS$ and $\sigma_\CS:\CS\tens_A\Omega^1\to \Omega^1\tens_A\CS$, and the above formulae now with $\omega\in \Omega^1$ replaced by $
\psi\in \CS$. We also need a {\em Clifford bimodule map} $\la:\Omega^1\tens_A\CS\to \CS$ which generalises the role of the gamma-matrices when constructing Dirac operators\cite{BegMa:spe}. In this case, we define $\dirac=\la\circ\nabla_\CS: \CS\to \CS$ as a minimal `geometric Dirac operator' in the formalism. We ideally would like $\doublenabla(\la)=0$ i.e., 
\begin{equation}\label{covla}  (\id\tens\la)\circ(\nabla\tens \id+ (\sigma\tens\id)(\id\tens\nabla_\CS))=\nabla_\CS\circ\la,\end{equation}
which says that $\la$ is covariantly constant in the sense that it intertwines the tensor product connection on $\Omega^1\tens_A\CS$ and the connection on $\CS$. \cite{BegMa:spe} also tentatively proposed a `Clifford relation' that  $\la$ extends to $\Omega^2\tens_A\CS\to \CS$ in such a way that 
\begin{equation}\label{cliffla} (\omega\wedge\eta)\la \psi:=\varphi(\omega\la(\eta\la \psi))-\kappa (\omega,\eta)\psi\end{equation}
for all $\omega,\eta\in \Omega^1$ and $\psi\in \CS$, and an optional bimodule automorphism $\varphi:\CS\to \CS$ and a nonzero constant $\kappa$. The motivation here is that this requirement, with $\varphi=\id, \kappa=1$, reduces to the Clifford algebra relations in the classical case. Both of these geometric conditions are part of the quantum Riemannian geometry and {\em not} part of the requirements for $\dirac$ to be part of a Connes spectral triple. However, (\ref{covla}) is relevant to the Laplacian on spinors defined as
\begin{equation}\label{lapS} \square_\CS= ((\ , \)\tens\id)\nabla_{\Omega^1\tens\CS}\nabla_\CS=((\ ,\ )\tens\id)(\nabla\tens \id+ (\sigma\tens\id)(\id\tens\nabla_\CS))\nabla_\CS\end{equation}
and when both (\ref{covla})-(\ref{cliffla}) hold, there is a Lichnerowicz theorem\cite[Prop.~8.45]{BegMa}  
\begin{equation}\label{lich} \varphi \circ\dirac^2=\kappa\square_\CS+\la\circ R_\CS, 
\end{equation}
where $R_{\CS}=(\extd\tens\id-\id\wedge\nabla_\CS)\nabla_\CS$ is the curvature of $\nabla_\CS$. 

For a Connes spectral triple, we instead need that $\CS$ has an antilinear skew-bimodule map $\CJ:\CS\to \CS$ in the role of `charge conjugation' and, if it exists, a bimodule map $\gamma:\CS\to \CS$ which we will refer to as an `even structure', both with certain properties needed to fit the algebraic part of Connes axioms\cite{Con} for free sign parameters $\eps,\eps',\eps''$. Certain patterns of signs have associated modulo 8 `KO dimension' $n$ and only even $n$ require $\gamma$, but in fact these restrictions come from classical motivations and we are not bound by them. Finally, we need a compatible sesquilinear inner product $\<\ ,\ \>$ on $\CS$ with respect to which we complete the latter into a Hilbert space as explained in \cite{BegMa:spe}\cite[Chap.~8.5]{BegMa}, where the inner product is denoted $(\!\!\<\ ,\ \>\!\!)$, such that  $\dirac$ is antihermitian, $\CJ$ an (antilinear) isometry and (if it exists) $\gamma$ hermitian. Here, $\imath\dirac$ is the hermitian Dirac operator in Connes axioms. Note that \cite{BegMa:spe}\cite[Chap.~8.5]{BegMa} uses a more careful notation in which antilinearity is expressed as linearity by the use of the conjugate bimodule $\bar{\CS}$ (defined as the same set as $\CS$ but conjugated actions of the field and of $A$) and $j(\psi)=\overline{\CJ(\psi)}$, etc. The  over-line here means to view the element of $\CS$ in $\bar {\CS}$. One should similarly distinguish between a sequilinear map and a linear map on $\bar\CS\tens\CS$. These distinctions are, however, only for book-keeping and we suppress them here.

\subsection{QRG realisation of Connes spectral triples for a central basis}\label{secspebasis}

We will be interested in the case where $A$ has trivial centre, $\Omega^1$ has a central basis $\{s^i\}$ and $\CS$ has a central basis $\{e^\alpha\}$.  Some general analysis for quantum Riemannian geometry on $\Omega^1$ in this parallelisable case is in \cite[Exercise~8.2 ]{BegMa} and was extended to construct QRG-realised spectral triples in \cite{LirMa2} but under more restrictive conditions than needed now. Over $\C$, we assume that $s^i{}^*=s^i$ are self-adjoint. Then `reality' of the metric translates to the matrix $g_{ij}$ of metric coefficient  in the basis being hermitian. If $\sigma$ is the flip map then the $*$-preserving condition on $\nabla$ translates to the Christoffel  symbols in the basis being real, otherwise the condition is more complicated. 

The further data for a Dirac operator starts with a  `Clifford action' $\la:\Omega^1\tens_A \CS\to \CS$  of the form
\begin{equation}\label{Cia} s^i\la e^\alpha= C^{i\alpha}{}_\beta e^\beta ,\quad  C^{i\alpha}{}_\beta\in\C\end{equation}
due to centrality of the bases and being a bimodule map.  We will also consider this as a collection of matrices $(C^i){}^\alpha{}_\beta=C^{i\alpha}{}_\beta$.  Next, we define the `charge conjugation' operator according to
\[ \CJ(a e^\alpha)= a^*J^\alpha{}_\beta e^\beta\]
\begin{equation}\label{JJ}J^\alpha{}_\beta{}^* J^\beta{}_\gamma=\eps\delta^\alpha{}_\gamma;\quad \overline{J}J=\eps\id,\quad \eps=\pm1,\quad J^\alpha{}_\beta\in \C,\end{equation}
where we also write the equation in a compact form for a matrix $J$. The over-line denotes complex conjugation of the entries.

We similarly need a bimodule connection which we write in the form
\[ \nabla_{\CS}e^\alpha=S^\alpha{}_{i\beta}s^i\tens e^\beta,\quad \sigma_\CS(e^\alpha\tens s^j)=\sigma_\CS{}^{\alpha j}{}_{i\beta }s^i\tens e^\beta,\quad \sigma_\CS{}^{\alpha j}{}_{i\beta}\in \C,\]
where now $S^\alpha{}_{i\beta}$ will have entries in $A$, not constants as in \cite{LirMa2} and $\sigma_{\CS}$ need not be the flip map (this is the same story as for $\Gamma,\sigma$ below). The abstract formulation of the `reality condition' on $\nabla_\CS$ in \cite{BegMa:spe} then becomes
\begin{equation}\label{SJ}
            S^\alpha{}_{i \beta}{}^* J^\beta{}_{\gamma}\sigma_\CS{}^{\gamma i}{}_{j\delta}=J^\alpha{}_\beta S^\beta{}_{j \delta} ;\quad \overline{S_j}J\sigma_\CS{}^j{}_i=J S_i.
         \end{equation}
Compatibility of $\CJ$ with the Clifford action in the general scheme of \cite{BegMa:spe} now becomes
         \begin{equation}\label{CJ}
            C^{i\alpha}{}_{\beta}{}^* J^\beta{}_{\delta} = \epsilon' J^\alpha{}_{\beta}\sigma_\CS{}^{\beta i}{}_{j\gamma}C^{j \gamma}{}_{\delta} ;\quad \eps'=\pm 1,\quad   \overline{C^i}J=\eps'J\sigma_\CS{}^i{}_j C^j.
          \end{equation} 

Finally, for an  `even structure' we also set $\gamma (e^\alpha) = \gamma^\alpha{}_{\beta} e^\beta$, here again with $\gamma^\alpha{}_{\beta} \in \C$ under our central basis and centre assumptions. Then  the properties needed in \cite{BegMa:spe} reduce to 
\begin{align}\label{gamsq}
          \gamma^\alpha{}_{\beta} \gamma^\beta{}_{\nu} &= \delta^{\alpha}{}_{\nu};\quad\quad\qquad \quad\gamma^2=\id,
         \\ \label{Cgam} 
          C^{i\alpha}{}_{\beta}\gamma^{\beta}{}_{\nu}&=-\gamma^\alpha{}_{\beta} C^{i\beta}{}_{\nu};\quad\qquad \{C^i,\gamma\}=0,          \\ \label{Jgam}
          \gamma^\alpha{}_{\beta}{}^* J^\beta{}_{\nu} &= \epsilon'' J^\alpha{}_{\beta} \gamma^\beta{}_{\nu};\quad\qquad \overline{\gamma}J=\eps''J \gamma,\quad\eps''=\pm1,
         \\ \label{Sgam}
          S^\alpha{}_{i\beta}\gamma^{\beta}{}_{\nu}&=\gamma^\alpha{}_{\beta} S^\beta{}_{i\nu};\quad\quad\quad [S_i,\gamma]=0.
        \end{align}

Here the sign is (\ref{Cgam}) is as for an even `KO-dimension'  $n$ but it turns out to be natural to adopt it for all cases. We note that the full set of conditions for $J$ is invariant under multiplication by a phase.  Finally, 
\begin{align*} \dirac(\psi_\alpha e^\alpha)&= \la(\nabla_\CS(\psi))=\la(\del_i \psi_\alpha s^i\tens e^\alpha+\psi_\beta\nabla_\CS e^\beta)=(\del_i\psi_\alpha+ \psi_\beta S^\beta{}_{i\alpha})C^{i\alpha}{}_\gamma e^\gamma\end{align*}
gives the resulting Dirac operator
\begin{equation}\label{dirac}   (\dirac\psi)_\gamma=(\del_i\psi_\alpha+ \psi_\beta S^\beta{}_{i\alpha})C^{i\alpha}{}_\gamma\end{equation}
in our conventions.

Beyond the above data for  a spectral triple at the local tensorial level, we have two further, quantum-geometric, requirements not present in Connes' formalism. The first is 
the covariance $\doublenabla(\la)=0$ condition:
\begin{equation} \label{dnablacliff} C^{i\alpha}{}_\beta S^\beta{}_{j\gamma}- \sigma^{ik}{}_{jl}  S^\alpha{}_{k\beta}C^{l\beta}{}_\gamma=-{1\over 2}\Gamma^i{}_{jk}C^{k\alpha}{}_\gamma;\quad C^iS_j-\sigma^{ik}{}_{jl}S_k C^l=-\tfrac{1}{2} \Gamma^{i}{}_{jk}C^k,\end{equation}
where 
\[ \nabla s^i=- {1\over 2}\Gamma^{i}_{jk}s^j\tens s^k,\quad \sigma(s^i\tens s^j)=\sigma^{ij}{}_{kl}s^k\tens s^l,\quad \sigma^{ij}{}_{kl}\in\C\]
in our conventions and under our assumptions. Unlike \cite{LirMa2}, we allow $\Gamma^i{}_{jk}\in A$ and $\sigma$ need not be the flip. The second quantum-geometric requirement is the compatibility of the Clifford action with $\Omega^2$, namely,
 \begin{equation}\label{cliffrel}  (s^i\wedge s^j)\la e^\alpha:=C^{j\alpha}{}_\beta C^{i\beta}{}_\gamma \varphi^\gamma{}_\delta e^\delta-\kappa g^{ij} e^\alpha \end{equation}
 is required to be well-defined. Here $g=g_{ij}s^i\tens s^j$ is the quantum metric with coefficients forming a hermitian matrix and $g^{ij}$ is its inverse, and $\varphi(e^\alpha)=\varphi^\alpha{}_\beta e^\beta$ for $\varphi^\alpha{}_\beta\in \C$ under our assumptions. Again unlike \cite{LirMa2}, we do not assume that the $s^i$ enjoy the usual Grassmann algebra, so this does not necessarily impose that the $C^i$ obey something resembling the usual Clifford algebra relations. We also  compute the spinor Laplacian
\begin{proposition}\label{propLapS} The spinor Laplacian in the case of the above bases is given by
\[ (\square_\CS \psi)_\alpha= \square \psi_\alpha+ \psi_\beta L{}^\beta{}_\alpha+ (\del_i\psi_\beta)(g^{ij}+g^{kl}\sigma^{ij}{}_{kl})S^\beta{}_{j\alpha},\]
where
\[\square a=(\ ,\ )\nabla((\del_j a)s^j)=g^{ij}\del_i\del_j a- {1\over 2} (\del_j a)g^{kl}\Gamma^j{}_{kl},\]
\[ L^\alpha{}_\beta=g^{ij}\del_i S^\alpha{}_{j\beta}-{1\over 2} S^\alpha{}_{j\beta} g^{kl} \Gamma^j{}_{kl}+g^{kl}\sigma^{ij}{}_{kl}S^\alpha{}_{i\gamma}S^\gamma{}_{j\beta},\]
for all $a\in A$ are the scalar Laplacian and coefficients defined by $\square_\CS(e^\alpha):=L^\alpha{}_\beta e^\beta$, respectively. 
\end{proposition}
\proof We use the derivation properties of the two connections and then their description in the basis,
\begin{align*}\square_\CS&(a e^\alpha)=(\ ,\ )_{12}\nabla_{\Omega^1\tens\CS}\nabla_\CS(a e^\alpha)=(\ ,\ )_{12}\nabla_{\Omega^1\tens\CS}((\del_j a) s^j\tens e^\alpha+a\nabla_\CS e^\alpha)\\
&= (\ ,\ )_{12}((\del_i\del_j a) s^i\tens s^j\tens e^\alpha+(\del_i a) s^i\tens S^\alpha{}_{j\beta}s^j\tens e^\beta+ (\del_j a)\nabla_{\Omega^1\tens \CS}(s^j\tens e^\alpha))+ a\square_\CS e^\alpha\\
&=g^{ij}(\del_i\del_j a) e^\alpha+g^{ij}(\del_i a)  S^\alpha{}_{j\beta} e^\beta+ (\del_j a)(\ ,\ )_{12}\nabla_{\Omega^1\tens \CS}(s^j\tens e^\alpha)+ a\square_\CS e^\alpha
\end{align*}
into which we put
\begin{align*} (\ ,\ )_{12}\nabla_{\Omega^1\tens \CS}(s^i\tens e^\alpha)&=(\ ,\ )(\nabla s^i\tens e^\alpha+ \sigma_{12}(s^i\tens S^\alpha{}_{j\beta}s^j \tens e^\beta)\\
&=-{1\over 2}\Gamma^i{}_{kl}g^{kl} e^\alpha+ \sigma^{ij}{}_{kl}g^{kl}S^\alpha{}_{j\beta}e^\beta.
\end{align*}
The first term of this combines with the first term of the first calculation to give $(\square a)e^\alpha$, and the other term goes into the expression stated for $\square_\CS$ on setting $a=\psi_\alpha$ and identifying coefficients in the answer. We also need
\begin{align*}\square_\CS(e^\alpha)&=(\ ,\ )_{12}\nabla_{\Omega^1\tens \CS}(S^\alpha{}_{i\beta}s^i\tens e^\beta)=(\del_j S^\alpha{}_{i\beta})g^{ji}e^\beta+ S^\alpha{}_{i\beta}(\ ,\ )_{12}\nabla_{\Omega^1\tens \CS}(s^i\tens e^\beta)
\end{align*}
into which we put the result of our second computation and then read off the coefficients $L^\alpha{}_\beta$.
\endproof
Next, if the Clifford relations condition (\ref{cliffrel}) is imposed then have an action of $\Omega^2$ on $\CS$ and hence $\la\circ R_\CS:\CS\to \CS$ makes sense. It is necessarily a left module map so it is enough to compute it on a basis element. 

\begin{lemma}\label{lemRS} Under the assumption of (\ref{cliffrel}) and $\nabla$ torsion free and with curvature coefficients defined by $(\la\circ R_\CS)(e^\alpha):=R_\CS{}^\alpha{}_\beta e^\beta$, we have
\[ R_\CS{}^\alpha{}_\eta =(\del_i S^\alpha{}_{j\beta}-{1\over 2}S^\alpha{}_{k\beta}\Gamma^k{}_{ij}- S^\alpha{}_{i\tau}S^\tau{}_{j\beta} )(C^{j\beta}{}_\gamma C^{i\gamma}{}_\delta \varphi^\delta{}_\eta -\kappa g^{ij} \delta^\beta{}_\eta).\]
\end{lemma}
\proof 
\begin{align*}(\la \circ R_\CS)(e^\alpha)&=\la (\extd\tens \id- \id\wedge\nabla_\CS)(S^\alpha{}_{i\beta} s^i\tens e^\beta)\\
&=(\del_i S^\alpha{}_{j\beta} s^i\wedge s^j+ S^\alpha{}_{k\beta} \extd s^k)\la e^\beta- \la(S^\alpha{}_{i\beta} s^i\wedge \nabla_\CS e^\beta)\\
&=(\del_i S^\alpha{}_{j\beta} s^i\wedge s^j+ S^\alpha{}_{k\beta} \extd s^k)\la e^\beta- S^\alpha{}_{i\tau}S^\tau{}_{j\beta} s^i\wedge s^j\la e^\beta
\end{align*}
into which we put torsion freeness $\extd s^k=\wedge \nabla s^k=-{1\over 2}\Gamma^k{}_{ij}s^i\wedge s^j$ and (\ref{cliffrel}), and then read off the coefficients as stated. \endproof

When we have a full quantum-geometric realisation so that both (\ref{dnablacliff})-(\ref{cliffrel}) hold then one has a Lichnerowicz formula\cite[Prop.~8.45]{BegMa}, which in our basis now reads
\begin{equation}\label{diracsqmat} (\dirac^2\psi)_\alpha \varphi^\alpha{}_\beta=\kappa (\square_\CS\psi)_\beta+ \psi_\alpha R_\CS{}^\alpha{}_\beta. \end{equation}
 This completes our analysis of the `local tensorial level' of the geometric realisation problem for a larger class of quantum geometries with central bases than in \cite{LirMa2}, as well as some of the related spinor quantum geometry. 

Next we note that the equation (\ref{JJ}) is a self-contained problem and can already be analysed in each dimension. When $\CS$ is 2-dimensional over the algebra, which will be the case of interest, it is easy enough to describe the possible $J$ explicitly up to a phase (which is not determined by the equations). 

\begin{lemma}\label{lemJ} Up to a phase in the 2D spinor bundle case, $J$ can be obtained with $r>0$ and  $z\in\C$ as either:
\[ {\rm (1):}\quad   J=\begin{pmatrix} z & r\\ {\eps-|z|^2\over r} & -\bar z\end{pmatrix},\quad {\rm (2):}\quad J=\begin{pmatrix} 1 & {\eps-1\over |z|^2}z \\ z &- {z\over\bar z}\end{pmatrix}\]
or its transpose. The $\eps=-1$ case of (2) needs $z\ne 0$ and up to a phase is also an instance to type (1). 
\end{lemma}
\proof Up to a phase we can assume the top right corner is real or zero, so 
  $J=\begin{pmatrix}a & r\\ c& d\end{pmatrix}$ or $J=\begin{pmatrix}a& 0\\ c& d\end{pmatrix}$ for complex parameters and a real $r>0$. Writing out (\ref{JJ}) we have in the first case looking off-diagonal in $\bar J J$ that $d=-\bar a$  and $(\bar c-c)a=0$. If $a\ne 0$, we need $c$ real and $c=(\eps-|a|^2)/r$ which gives the form shown. If $a=0$ then we are led to $c=\eps/r$.  For the other case, from the diagonal of $\bar JJ$ we have $|a|^2=|d|^2=\eps$ hence $\eps=1$ and, up to a phase, we can assume $a=1$ also.  In this case off diagonal we need that $d=-c/\bar c$ as shown.  We note, however, that those cases of type (1) for general $r$ and $\eps=-1$ that have an element of modulus 1 in the top left are equal up to a phase to the case (2) with  $\eps=-1$, so this is redundant if we want an exclusive classification. Also for $\eps=1$ in type (2), we can send $z\to 0$ along a ray, which up to a phase then overlaps with $z$ a phase and $r=0$ as a limit of type (1). \endproof
  
 We also note that our (\ref{JJ})--(\ref{CJ}) are invariant under global conjugation by an invertible $u$, 
 \begin{equation}\label{unitaryequiv} C^i\mapsto uC^i u^{-1},\quad S_i\mapsto uS_i u^{-1},\quad \sigma_\CS{}^i{}_j\mapsto u\sigma_\CS{}^i{}_j u^{-1},\quad J\mapsto \overline{u}Ju^{-1},\end{equation}
 while $\sigma^{ij}{}_{kl}\in \C$ and hence is invariant. If we allow this then we can conjugate the type (1) $J$ by a diagonal matrix of determinant 1 to set $r=1$. Hence, in calculations for the type (1) case, we will generally only consider $r=1$ for simplicity. Then for both types, $J$ is given by a single complex parameter $z$. 
 
 It remains to address the Hilbert space side of a spectral triple. For the inner product, we will assume  a reference positive linear functional $\int: A\to \C$ and set 
 \begin{equation}\label{hilbS}\< \phi_\alpha e^\alpha, \psi_\beta e^\beta\>=\int \phi_\alpha^*\mu^{\alpha\beta}\psi_\beta\end{equation}
 for some positive hermitian matrix $\mu$ as `measure' and $\phi,\psi\in\CS$. We use this to complete $\CS$ to a Hilbert space. One then needs to check that $\CJ$ is an antilinear isometry, which, adapting the calculation in \cite{LirMa2} and if we assume that $\int$ is a trace, amounts to
 \begin{equation}\label{muJ} \eps \mu^{\alpha\gamma}J^{\beta}{}_\gamma=J^\alpha{}_\gamma \mu^{\gamma\beta}.\end{equation}
The simplest case, which will also suffice for our purposes, is the identity matrix $\mu^{\alpha\beta}=\delta^{\alpha\beta}$ as in \cite{LirMa2}. Then the condition on $J$ from this is that it be symmetric if $\eps=1$ and antisymmetric if $\eps=-1$.    For a spectral triple, we  need that $\dirac$ from (\ref{dirac}) and the operator $\gamma$ (if present) are antihermitian and hermitian respectively, which it is easiest at this point to check directly.

\subsection{Case of inner calculus.} \label{secspeinn} We now specialise to the case where the calculus $\Omega^1$ is inner with $\theta=\theta_i s^i$ for $\theta_i\in A$. Then by \cite[Prop.~8.11]{BegMa}, we have $\nabla s^i=\theta\tens s^i-\sigma(s^i\tens\theta)-\alpha(s^i)$ for some bimodule map $\alpha(s^i)=\alpha^i{}_{jk}s^j\tens s^k$, which translates to
\begin{equation}\label{Gtheta} -{1\over 2}\Gamma^i{}_{jk}=\theta_j\delta^i_k -\theta_m\sigma^{im}{}_{jk}+\alpha^i{}_{jk};\quad \alpha^i{}_{jk}\in \C.\end{equation}
Similarly for $\nabla_\CS$ (see \cite[Sec.~6.4]{MaSim}), we will have $\nabla_\CS=\theta\tens(\ )- \sigma_\CS((\ )\tens\theta)+\alpha_\CS$ for  $\sigma_\CS$ as above and another bimodule map,
\begin{equation}\label{SthetaA} \alpha_\CS(e^\alpha)=A^\alpha{}_{i\beta}s^i\tens e^\beta,\quad  A^\alpha{}_{i\beta}\in \C,\quad S^\alpha{}_{i\beta}=\theta_i\delta^\alpha_\beta-\theta_j\sigma_\CS{}^{\alpha j}{}_{i\beta}+ A^\alpha{}_{i\beta}.\end{equation}
For an inner calculus, one says that $\nabla,\nabla_\CS$ are inner if they are given entirely by $\sigma,\sigma_\CS$ respectively, i.e. $\alpha^i{}_{jk}=0$ and $A^\alpha{}_{i\beta}=0$ respectively.  

 \begin{lemma}\label{leminn} Let $\nabla$ be inner and write $(\sigma_\CS{}^i{}_j)^\alpha{}_\beta:=\sigma_\CS{}^{\alpha i}{}_{j\beta}$ for the matrix versions.  (1) If $\nabla_\CS$ is inner then equation (\ref{dnablacliff}) holds if and only if
\[C^{i\alpha}{}_\beta \sigma_\CS{}^{\beta k}{}_{j\gamma}=\sigma^{im}{}_{jl}  \sigma_\CS{}^{\alpha k}{}_{m\beta}C^{l\beta}{}_\gamma;\quad C^i\sigma_\CS{}^k{}_j=\sigma^{im}{}_{jl}\sigma_\CS{}^k{}_m C^l\]
and equation  (\ref{SJ}) holds if and only if
\[  J^\alpha{}_\delta \delta_{jk}=\overline{\sigma_\CS{}^{\alpha k}{}_{i\beta}}J^\beta{}_\gamma \sigma_\CS{}^{\gamma i}{}_{j\delta};\quad  \delta_{jk}J=\overline{\sigma_\CS{}^k{}_i}J \sigma_\CS{}^i{}_j.\]
(2) If $\nabla_\CS$ is general with $A^\alpha{}_{i\beta}$ then (\ref{dnablacliff})-(\ref{SJ}) are equivalent to the conditions in part (1) and additionally
\[ C^iA_j=\sigma^{ik}{}_{jl}A_kC^l,\quad J A_i=\overline{A_i}J\sigma_\CS{}^i{}_j.\]
 \end{lemma}
\proof (1) We take $A_i=0$ and the inner form of $\Gamma^i{}_{jk}$ and substitute into (\ref{dnablacliff}) to obtain
\[ \theta_m C^{i\alpha}{}_\beta \sigma_\CS{}^{\beta m}{}_{j\gamma}=\theta_m\sigma^{ik}{}_{jl}  \sigma_\CS{}^{\alpha m}{}_{k\beta}C^{l\beta}{}_\gamma.\]
Now note that by the surjectivity assumption on a calculus, for each $i$ there exist $a^i_\mu , b^i_\mu\in A$ such that $\sum_\mu a^i_\mu[\theta_j,b^i_\mu]s^j=s^i$. Since the $\{s^i\}$ are a basis, this means $\sum_\mu a^i_\mu[\theta_j,b^i_\mu]=\delta^i_j$. Hence, if $
\sum_j c_j\theta_j=0$ in $A$ for some $c_i\in\C$ then $c_j=0$, so the $\theta_j$ are linearly independent whenever we have a central basis. Hence, we can cancel the $\theta_m$ to obtain the result. Similarly, (\ref{SJ}) becomes the condition stated on putting in $A=0$ and noting that $\theta^*=-\theta$ so $\theta_j^*=-\theta_j$ given our central self-adjoint basis. (2) For the last part, we repeat the above with the extra $A$ terms on both sides. These are valued in multiples of $1$ and this can be taken as linearly independent of the $\theta_i$ since, given the central basis, any constant parts of $\theta_i$ would not affect commutators in computing $\extd f$. Hence, these parts of the equations decouple and appear as shown in matrix form. 
\endproof

For an even spectral triple, the condition (\ref{Sgam}) similarly translates in the inner case to
\begin{equation}\label{Sgaminn}  [\sigma_\CS{}^i{}_j,\gamma]=0,\quad [A_i,\gamma]=0\end{equation}
again using that the $\theta_i\in A$ are linearly independent and can be taken as linearly independent of 1.

\section{Dirac operator on noncommutative torus $\C_\theta[\T^2]$} \label{sector}

The well-known noncommutative torus at an algebraic level is the $*$-algebra $\C_{\theta}[\T^2]$ generated by unitary $u,v$ with relations $v u = e^{\imath \theta}u v.$ Here $u^*=u^{-1}$ and $v^*=v^{-1}$ and $\theta$ a real parameter (theres is no confusion with inner structures as
the calculus will not be inner). The Dirac operator here was one of the first constructions on the noncommutative torus, but here we revisit with the question of which quantum Riemannian geometries can realise it. We assume that $\theta$ is irrational, in which case the algebra has trivial centre. 
  
\subsection{Recap of WQRG of $\C_{\theta}[\T^2]$}\label{secqgtorus} 

  For $\Omega^1(\C_\theta[\T^2])$, the standard differential calculus has
  \[( \extd u )u=u\extd u,\quad (\extd v)v=v\extd v,\quad (\extd v)u=e^{\imath\theta}u\extd v,\quad (\extd u)v=e^{-\imath\theta}v\extd u\]
  with maximal prolongation
   \begin{equation*}
        \extd u \wedge \extd u =0, \quad \quad \extd v \wedge \extd v=0, \quad \quad \extd v \wedge \extd u = -e^{\imath\theta} \extd u \wedge \extd v
    \end{equation*}
  as explained e.g. in \cite[Chap. 1]{BegMa}. The first thing we do is choose a self-adjoint basis 
  \[ s^1=-\imath u^{-1}\extd u,\quad s^2=-\imath v^{-1}\extd v.\]
   Then one can check that  these are central and Grassmann with $\extd s^i=0$. The partial derivatives defined by $\extd f=(\del_i f)s^i$ for all $f\in\C_{\theta}[\T^2]$  are the standard derivations
   \[ \del_1(u^mv^n)=\imath m u^m v^n,\quad \del_2(u^m v^n)=\imath n u^m v^n.\]

    Next, a (generalised, i.e. not-necessarily quantum symmetric) quantum metric has the form    \begin{equation*}
        g=c_1 s^1 \otimes s^1 + c_2 s^2 \otimes s^2 + c_3 s^1\otimes s^2 + c_4 s^2 \otimes s^1
    \end{equation*}
   for $c_1,c_2,c_3,c_4 \in \C$, since the metric has to be central to be bimodule-invertible. We suppose non-degeneracy, which amounts to $c_1c_2\ne c_3 c_4$, and focus on the quantum symmetric case $c_4=c_3$. In this case the `reality' of the metric
    says that $c_i$ are real. 
    
  Similarly, it is shown \cite[Example~8.16]{BegMa} that a torsion free bimodule connection must have the form
  \begin{equation}\label{torfree1}
        \nabla s^1 = h^1{}_1 s^1\otimes s^1 + h^1{}_2 s^2\otimes s^2 + h^1{}_3 (s^1\otimes s^2 +s^2 \otimes s^1 ) ,
    \end{equation}
    \begin{equation}\label{torfree2}  
        \nabla s^2= h^2{}_1 s^1\otimes s^1 + h^2{}_2 s^2\otimes s^2 + h^2{}_3 (s^1\otimes s^2 +s^2 \otimes s^1 )
    \end{equation}
    for  constants $h^i{}_j\in \R$ and $\sigma={\rm flip}$ on the basis 1-forms. Putting this form in the requirements for a QLC, we find that $\nabla s^i=0$ is the only solution in the quantum symmetric case. More generally, 
    we can ask only for the torsion and cotorsion to vanish, which happens in the quantum symmetric case when 
 \begin{equation}\label{torwqlc} c_2 h^2{}_1 = c_1 h^1{}_3 + c_3 (h^2{}_3 -h^1{}_1),\quad  c_2 h^2{}_3 = c_1 h^1{}_2 +c_3 (h^2{}_2 -h^1{}_3 ),\end{equation}
  giving a 4-parameter moduli of WQLC's\cite{BegMa}. These have curvature
      \begin{equation*}
        R_{\nabla} (s^1)= \rho^1{}_j s^1 \wedge s^2 \otimes s^j, \quad R_{\nabla} (s^2)= \rho^2{}_j s^1 \wedge s^2 \otimes s^j,
    \end{equation*}  
    where\cite{BegMa}
    \begin{equation*}
        \rho = S \begin{pmatrix} c_3 & c_2 \\ -c_1 & -c_3 \end{pmatrix},  \quad  S=\begin{cases} \frac{h^1{}_2 h^2{}_1 - h^1{}_3 h^2{}_3}{c_3}&{\rm if\ }c_3\ne 0\\ {h^2{}_3(h^2{}_3-h^1{}_1)\over c_1}+ {h^1{}_3(h^1{}_3-h^2{}_2)\over c_2}&{\rm if\ }c_3=0.\end{cases}
    \end{equation*}
 Here $S$ is the Ricci scalar and the ${\rm Ricci}={S g\over 2}$ so that the `obvious' Einstein tensor vanishes. We used the usual antisymmetric lift of $s^1\wedge s^1$. We also will need a spinor connection:
 
 \begin{lemma} For a bimodule connection on $\CS$ with central  basis $\{e^\alpha\}$ over $\C_\theta[\T^2]$, $\sigma_\S$ has to be the flip map and $S_i^\alpha{}_\beta\in \C$.
 \end{lemma}
 \proof Computing $\nabla_\CS(e^\alpha a)=\nabla_\CS(ae^\alpha)$ using the two Leibniz rules, we have that
\[ (\sigma_\CS{}^{\alpha j}{}_{i\beta}-\delta^j{}_i\delta^\alpha{}_\beta)\del_j a=[a, S_i^\alpha{}_\beta]\]
for all $a=u^m v^n$, say. But $\C_\theta[\T^2]$ has a $\Z\times\Z$ grading according to the powers of $u,v$ and the left hand side has the same grade $(m,n)$ as $a$. If we write $S_i^\alpha{}_\beta=\sum_{a,b} S^{(a,b)}_i{}^\alpha{}_\beta u^a v^b$ then the $(a,b)$ term has grade $(m+a,n+b)$ on the right hand side of the displayed equation. Hence, we conclude that $S^{(a,b)}_i$ can only contribute to the right hand side if $(a,b)=(0,0)$ and in this case the right hand side is zero. Hence, $\sigma_\CS$ is necessarily the flip map. It then follows that $S_i^\alpha{}_\beta\in \C$ as the centre of the algebra is trivial. \endproof

This also requires that $\nabla$ has to have $\sigma$ the flip map and $\C$-number Christoffel symbols as in \cite{BegMa}, but now argued more generally.

\subsection{Uniqueness of geometric spectral triple on $\C_\theta[\T^2]$} 

To be concrete, we take $c_1=c_2=1$ and $c_3=0$ so that $g=s^1\tens s^1+s^2\tens s^2$. In this case, given the Grassmann form of $\wedge$,  the Clifford relations (\ref{cliffrel}) with $\varphi=\id$ and $\kappa=1$ become $C^iC^j+C^jC^i=2g^{ij}$, which we solve for  2-dimensional $\CS$ and up to unitary equivalence by the usual choice
\[ C^i=\sigma^i;\quad  \sigma^1=\begin{pmatrix}0&1\\1&0\end{pmatrix},\quad \sigma^2=\begin{pmatrix}0&-\imath\\ 
\imath&0\end{pmatrix},\quad \sigma^3=\begin{pmatrix}1&0\\0&-1\end{pmatrix},\]
where we recalled the Pauli matrices, needed here for $i=1,2$ 

\begin{lemma} For $C^i=\sigma^i$, $J$ has the type (1) form
\[ \eps=\eps'=1:\quad J=q\sigma^1,\quad  \eps=\eps'=-1:\quad J=q\sigma^2\]
with $|q|=1$.  In either case, we also have $\gamma=\sigma^3$  with $\eps''=-1$, thereby solving all our equations at the local tensorial level that do not directly involve $S_i$. 
\end{lemma} 
\proof Given that $\sigma_\CS$ has to be the flip map, (\ref{CJ}) reduces to $\overline{C^i} J=\eps' J C^i$ 
as in \cite{LirMa2}. Also  (\ref{JJ}) has the standard form $\bar JJ=\eps\id$ with solutions in Lemma~\ref{lemJ}. In fact, the 
result is also easily found directly by writing $J=j_0\id+j_i\sigma^i$ and noting that $\overline{\sigma^2}=-\sigma^2$ while the other $\sigma^i$ have real entries. So $\sigma^1J=\eps' J\sigma^1$ and $\sigma^2J=-\eps'J\sigma^2$. We use the Pauli matrix algebra to reduce to the basis and thereby  deduce that if $\eps'=1$ then only $j_1$ is allowed, and if $\eps'=-1$ then only $j_2$ is allowed. We then check the other condition and find that $\eps=\eps'$ and $|j_1|=1$ or $|j_2|=1$ respectively. $J$ is always determined only up to a phase, which we denote $q$. It is then easy to solve (\ref{gamsq})-(\ref{Cgam}) to find $\gamma=\sigma^3$ as anticommuting with $\sigma^1,\sigma^2$. In this case (\ref{Jgam}) also holds provided $\eps''=-1$, for either form of $J$. \endproof

We next have to solve for $\nabla_\CS$ to be compatible with $\nabla, C^i, J$ so as to  fully solve the  requirements for a geometric spectral triple at the local tensorial level. 

\begin{lemma}\label{torSsolve} For $C^i=\sigma^i$, the requirements of a geometrically realised spectral triple require the WQLC to be the QLC $\nabla s^i=0$ and $\nabla_\CS e^\alpha = d_i s^i\tens e^\alpha$ for some real coefficients $d_i$. We also have  $\gamma$ at this level  as in the preceding lemma.  
\end{lemma}
\proof  Since  $\sigma_\CS$ is the flip, (\ref{SJ}) is  $\bar S_i J= J S_i$. 
Next, since the braiding $\sigma$ of the connection is also trivial, we need $[C^i,S_j]=-{1\over 2}\Gamma^i{}_{jk}C^k$. But writing $\Gamma$ in terms of $h^i{}_j$ in this case (where we solve the WQLC conditions (\ref{torwqlc}) by $h^1{}_3=h^2{}_1$ and $h^2{}_3=h^1{}_2$), we have
\[ -{1\over 2}\Gamma^i{}_{jk}=h^i{}_j\delta_{jk}+ h^{\bar i}{}_i \delta_{\bar j,k}\]
where $\bar 1=2$ and $\bar 2=1$ and $i,j,k$ are now restricted to $1,2$. Hence the covariance condition (\ref{dnablacliff}) is 
\[ [C^i,S_j]=h^i{}_j C^j+ h^{\bar i}{}_i C^{\bar j}.\]
Using Mathematica to search among this 4-parameter moduli of $\nabla$, we find that there are no solutions unless $h^i{}_j=0$, i.e. the trivial QLC, and in this case $S_1$ and $S_2$ have to be real multiples of the
identity since they have to commute with both $C^i$ and commute up to complex conjugation with $J$. In this case, the remaining (\ref{Sgam}) is automatic. \endproof

Hence the only geometrically realised Dirac operator at the local tensorial level for the standard Euclidean metric and $\nabla$ a WQLC are
\begin{equation}\label{tordirac} \dirac(\psi_\alpha e^\alpha)=(\del_i\psi_\alpha s^i)\la e^\alpha+ \psi_\alpha d_i s^i\la e^\alpha =\sigma^i{}^\alpha{}_\beta ((\del_i+ d_i)\psi_\alpha) e^\beta. \end{equation}
For an actual spectral triple, we still need to construct the Hilbert space, which we do with respect to the state $\int u^mv^n=\delta_{m,0}\delta_{n,0}$  (this would be integration on the classical torus). We use this to complete $\CS=\C_\theta[\T^2]\tens\C^2$ to a Hilbert space via (\ref{hilbS}) with $\mu=\id$ the identity matrix.  As $J$ is symmetric if $\eps=1$ and antisymmetric if $\eps=-1$, we already know that $\CJ$ is an antilinear isometry in either case for this Hilbert space. The remaining arguments are similar to \cite{LirMa2}.

\begin{theorem}\label{thmtorspec} Up to unitary equivalence, $(\dirac\psi)_\beta=(\del_i\psi_\alpha)\sigma^i{}^\alpha{}_\beta$ is the only possibility for a geometrically realised spectral triple on the noncommutative torus for the Euclidean metric,  a WQLC and the standard Hilbert space structure on $\CS$. Moreover, we are forced to $\nabla=\nabla_\CS=0$ on the bases. Here $\dirac$ forms a spectral triple of $KO$-dimension  $n=2$ if we choose $\eps=\eps'=\eps''=-1$ and $n=1$ if $\eps=\eps'=1$ ignoring $\gamma$. 
\end{theorem}
\proof  We have done the work in the two lemmas and just need to check if (\ref{tordirac}) is antihermitian and $\gamma$ is hermitian. For the former, 
\[ \<\dirac \psi,\phi\>+\<\psi,\dirac\phi\>= \int (\sigma^i{}^\alpha{}_\beta ((\del_i+ d_i)\psi_\alpha)^*\phi_\beta+ \psi_\alpha^* \sigma^i{}^\beta{}_\alpha (\del_i+ d_i)\phi_\beta =2\int d_i\sigma^i{}^\beta{}_\alpha \psi^*_\alpha\phi_\beta.\]
The terms without $d_i$ vanish since $\sigma^i$ are hermitian, $\del_i$ are derivations and $\int\del_i(\psi_\alpha^*\phi_\beta)=0$ because
$\int\del_i(u^mv^n)=n\int u^m v^n=0$ for all $m,n$. For our expression to be zero for all $\phi,\psi$, this requires $d_i=0$. For $\gamma$, we have more easily $\<\gamma\psi,\phi\>=\int (\sigma^3{}^\alpha{}_\beta\psi_\alpha)^*\phi_\beta=\<\psi,\gamma\phi\>=0$  because $\sigma^3$ is hermitian. Finally, looking in the period 8 table of signs, and allowing that $\dirac$ is antihermitian (the actual operator in the spectral triple is $\imath\dirac$) which flips the required value of $\eps'$, we see that we can meet the conditions for an even spectral triple with $n=2$ if we choose $\eps=\eps'=-1$. In the other case, we fit with $n=1$ and can choose to ignore $\gamma$. \endproof

We have, in particular, recovered the standard $n=2$ spectral triple on the noncommutative torus, including a suitable charge conjugation $J$. Although not a new spectral triple, we see how it is forced from the quantum geometry up to unitarity for the given choice of sign (ditto for the other sign). Other symmetric metrics will be equivalent after an orthogonal transformation to the Euclidean one and can therefore be handled in a manner similar to our treatment for the fuzzy sphere in  \cite{LirMa2}, again giving a unique answer up to unitary equivalence for a given choice of  signs.  Asymmetric metrics, however, won't be compatible with (\ref{cliffrel}) since, due to the Grassmann algebra form of the wedge product, we have to solve $C^iC^j+C^jC^i=2g^{ij}$, which requires $g^{ij}$ to be symmetric. 

We also have
\[ (\dirac^2\psi)_\alpha=\sum_{i=1}^2\del_i^2\psi_\alpha=\square \psi_\alpha\]
for the standard Dirac operator in Theorem~\ref{thmtorspec} as an easy instance of the general Lichnerowicz formula, which applies since both (\ref{covla})-(\ref{cliffla}) hold, in the present case with $\nabla_\CS$ flat.

\subsection{Similar results for more general $\nabla$}\label{sectorext}

Sticking with our Euclidean metric, which enters into the choice of  $C^i=\sigma^i$ for a 2-dimensional $\CS$, we can let $\nabla$ be any torsion free bimodule connection (\ref{torfree1})-(\ref{torfree2}), which is a 6-parameter moduli of connections. 

\begin{proposition}\label{propext1} Extending the geometric realisation to allow any torsion free $\nabla$ leads to the same answer as in Theorem~\ref{thmtorspec} with $\nabla=\nabla_\CS=0$ on the bases. Allowing $\nabla$ any bimodule connection leads to a real 4-parameter moduli of $\nabla,\nabla_\CS$ but $\dirac$ itself of the same form as in Lemma~\ref{torSsolve}, resulting in the same spectral triple(s) as in Theorem~\ref{thmtorspec}.\end{proposition}
\proof  If we drop the WQLC condition but keep the torsion free form and repeat the calculation in Lemma~\ref{torSsolve} then we do not obtain any more solutions, we are again forced to $\nabla=0$ in the basis. If, however, we allow all 8 values of $\Gamma^i{}_{jk}$  for a general bimodule connection, then we find solutions to (\ref{dnablacliff}) provided $\nabla$ has the 2-parameter form
\[  \nabla s^1=      h^1{}_2 s^2\otimes s^2 -h^2{}_1s^1\otimes s^2  ,\quad     \nabla s^2= h^2{}_1 s^1\otimes s^1 - h^1{}_2 s^2 \otimes s^1, \]
which we see is not torsion free unless it is zero, and 
\[ S_1=\begin{pmatrix}a_1 & 0 \\ 0 & a_1+\imath h^2{}_1\end{pmatrix},\quad S_2=\begin{pmatrix}a_2 & 0 \\ 0 & a_2-\imath h^1{}_2\end{pmatrix}\]
for free parameters $a_i\in \C$.  Then (\ref{SJ}) for either choice of $J$ requires the reality condition $\bar a_1=a_1+\imath h^2{}_1$ and $\bar a_2=a_2-
\imath h^1{}_2$, hence two further real parameters $d_i\in R$ with  
\[ a_1= d_1 - {\imath\over 2} h^2{}_1,\quad a_2=d_2+{\imath\over 2}h^1{}_2;\quad  S_1=d_1\id-{\imath\over 2}h^2{}_1\sigma^3,\quad S_2=d_2\id+{\imath \over 2}h^1{}_2\sigma^3\]
as the $h^i{}_j$ are real. One can check that (\ref{Sgam}) still holds, so we have the algebraic properties for an even spectral triple, although only relevant for one of the sign choices.

Finally, we do a similar analysis as before for $\dirac$. From (\ref{dirac}), we have
\begin{align*}  (\dirac\psi)_\gamma&=(\del_i\psi_\alpha+ d_i\psi_\alpha)\sigma^{i\alpha}{}_\gamma-{\imath \over 2}h^2{}_1\psi_\beta(\sigma^3\sigma^1)^\beta{}_\gamma  +{\imath \over 2}h^1{}_2  \psi_\alpha(\sigma^3\sigma^2)^\alpha{}_\gamma\\
&= (\del_i\psi_\alpha+ d_i\psi_\alpha)\sigma^{i\alpha}{}_\gamma+ {\psi_\alpha\over 2}(h^2{}_1\sigma^2+h^1{}_2\sigma^1)^\alpha{}_\gamma= (\del_i\psi_\alpha+ d'_i\psi_\alpha)\sigma^{i\alpha}{}_\gamma\end{align*}
where $d_i'=d_i+h^i{}_{\bar i}/2$. This is of the same form as in Lemma~\ref{torSsolve} so is not antihermitian unless $d'_i=0$. The $\gamma$ is the same as before and hermitian. \endproof

We can go further and  drop the $\doublenabla(\la)=0$ covariance condition (\ref{dnablacliff}) entirely. Then $\nabla$ does not enter the construction at all and we can search among all bimodule connections $\nabla_\CS$. This extends the results of the preceding proposition to an even bigger moduli. 

\begin{proposition} Geometric realisation at the local tensorial level without imposing $\doublenabla(\la)=0$ leads to an real 8-parameter moduli of $\nabla_\CS$, with a real 4-parameter submoduli that are compatible with our chosen $\gamma$. The latter results in the same $\dirac$ as in   Lemma~\ref{torSsolve} and the same spectral triple(s) as in Theorem~\ref{thmtorspec}, while more generally we obtain a 1-parameter moduli of spectral triples with real parameter $m$, 
\[ \eps=\eps'=1: \quad \dirac \psi=\del_i\psi\sigma^i+\imath m \psi \sigma^3;\quad 
\eps=\eps'=-1,\quad \dirac\psi =\del_i\psi\sigma^i+\imath m \psi.\]
\end{proposition} 
\proof If we search among all bimodule connections then we only have to solve  the reality condition (\ref{SJ}), which  tells us that each of  $S_1,S_2$ are of the form
\[\eps=\eps'=1:\quad  S_i=\begin{pmatrix} a_i& b_i\\ \bar b_i & \bar a_i\end{pmatrix};\quad \eps=\eps'=-1:\quad S_i=\begin{pmatrix} a_i& b_i\\ -\bar b_i & \bar a_i\end{pmatrix}\]
 for the two choices of $J$. This is a real 8-parameter moduli of spinor connections.  In either case, the condition (\ref{Sgam}) for an for compatibility with $\gamma=\sigma^3$ holds iff $b_i=0$, which cuts us down to 4 real dimensions. We also write our general $S_i$ in terms of Pauli matrices and combine with the $C^i$ as in the proof of Proposition~\ref{propext1}. The only contributions to $\dirac$ that are antihermitian are of the form stated with $m={\rm Im}(b_1)+{\rm Re}(b_2)$, which is, however, zero if we want an even spectral triple.  \endproof

The two new spectral triples here could be viewed as $n=1$ and $n=3$ respectively according to the signs (after flipping $\eps'$). This proposition represents the most general spectral triple of geometric type on the noncommutative torus under the assumption of the standard Clifford structure for a spinor space of two dimensions.

\section{Dirac operator on matrices $M_2(\C)$} \label{secM2}

We let $\Omega^1(M_2(\C))$ be the standard 2D calculus as in \cite{BegMa:spe}\cite[Chapter~1]{BegMa} with $s,t$ central and $s^*=-t$. We first convert this to a more physical self-adjoint basis
\[ s^1={(s+t)\over\imath},\quad  s^2=s-t;\quad s={1\over 2}(\imath s^1+s^2),\quad t={1\over 2}(\imath s^1-s^2),\]
so that $s^i{}^*=s^i$. The extension to higher forms is by $s^2=t^2=0$, $s\wedge t=t\wedge s$ with top form ${\rm Vol}:=2\imath s\wedge t$ which is self-adjoint in the $*$-calculus. In terms of the $s^i$, the exterior algebra, which is not Grassmann, is 
\begin{equation}\label{M2wedge} s^1\wedge s^1=s^2\wedge s^2=\imath{\rm Vol},\quad s^1\wedge s^2=s^2\wedge s^1=0.\end{equation}
The exterior algebra is inner with $\extd=[\theta,\ \}$ (graded-commutator) which in our basis converts to 
\begin{equation}\label{M2theta} \theta=E_{12}s+E_{21}t={\imath\over 2}(\sigma^1s^1+\sigma^2s^2),\end{equation}
where we use two of the Pauli matrices. Using all three Pauli matrices and $1=\id$ as a basis of $M_2(\C)$, the differentials are then
\begin{equation}\label{M2extd}  \extd \sigma^1=\sigma^3s^2,\quad \extd\sigma^2=-\sigma^3 s^1,\quad \extd\sigma^3=\sigma^2 s^1-\sigma^1s^2,\quad \extd s^i=-\sigma^i{\rm Vol} \end{equation}
as equivalent in our $s^i$ terms to the formulae in \cite{BegMa} in terms of $s,t$.  The corresponding partial derivatives defined by $\extd f=(\del_i f)s^i$ are 
\begin{equation}\label{M2deli} \del_i\sigma^j=-\eps_{ijk} \sigma^k,\quad \del_i={\imath\over 2}[\sigma^i,\ ]\end{equation}
for $i=1,2$ and $j,k=1,2,3$.  We will use this differential calculus for both the QRG in this section and the alternate QRG in the next section. In both cases, $\Gamma^i{}_{jk}$ is {\em not constant} even though the metric has constant coefficients, which makes the analysis fundamentally harder than the fuzzy sphere case in \cite{LirMa2} or the noncommutative torus case in Section~\ref{sector}.

 \subsection{Standard quantum Riemannian geometry on $M_2(\C)$}\label{secM2qrg}
 
 The standard quantum metric $M_2(\C)$ with its 2D calculus in  \cite[Example 8.13]{BegMa} (times 2) is $g=2(s\tens t-t\tens s)$, which in our self-adjoint basis becomes
 \begin{equation}\label{M2g}  g=\imath(s^2\tens s^1-s^1\tens s^2).\end{equation}
 The QLC is not unique, but there is a known 4-parameter moduli of $*$-preserving QLCs. Within this, we focus for simplicity on 1-parameter sub-moduli space
 \[ \nabla s=2\theta\tens s + \mu E_{21}(s\tens t-t\tens s) + \mu E_{12}(t\tens t-s\tens s)\]
 \[ \nabla t=2\theta\tens t + \mu E_{12}(s\tens t-t\tens s) + \mu E_{21}(t\tens t-s\tens s)\]
 \[ \sigma(s\tens s)=(\mu-1)s\tens s-\mu t\tens t,\quad \sigma(s\tens t)=-\mu s\tens t+ (\mu-1)t\tens s\]
 \[ \sigma(t\tens s)=(-\mu-1)s\tens t+\mu t\tens s,\quad \sigma(t\tens t)=\mu s\tens s+(-\mu-1)t\tens t\]
 with one real parameter $\mu$. This is given by setting $\nu=-\mu$ and $\beta=-\alpha\to 0$ in \cite[Example 8.13]{BegMa}.  In terms of $s^i$, the connection appears as 
 \[ \nabla s^i=2\theta\tens s^i-\mu \eps_{ik} \sigma^i s^i\tens s^k,\quad \sigma(s^i\tens s^j)= - s^j\tens s^i -\delta_{ij}2 \imath \mu \eps_{ik}s^i\tens s^k\]
 (sum over $k$). We also compute the curvature:

 \begin{proposition} The curvature of the above 1-parameter QLC for the metric (\ref{M2g}) is
 \[ R_\nabla(s^i)=-\mu\eps_{ij}{\rm Vol}\tens s^j   \]
 with Ricci curvature and Ricci scalar
 \[  {\rm Ricci}={ \mu\over 2} g,\quad S=-\mu.\]
 \end{proposition}
 \proof For the Riemann curvature we use that the calculus is inner to obtain
 \begin{align*} R_\nabla(s^i)&=(\extd\tens\id-\id\wedge\nabla)(2\theta \tens s^i-\mu \eps_{ik}\sigma^i s^i\tens s^k)\\
 &=2\theta^2\tens s^i-2\theta(2\theta \tens s^i-\mu \eps_{ik}\sigma^i s^i\tens s^k)- \mu\eps_{ik}\extd(\sigma^i s^i)\tens s^k\\
 &\quad +\mu \eps_{ik}\sigma^is^i(2\theta \tens s^k-\mu \eps_{km}\sigma^k s^k\tens s^m)\\
 &=\mu\eps_{ik}\extd(\sigma^i s^i)\tens s^k-\mu^2 \eps_{ik}\eps_{km}\sigma^i\sigma^k s^i s^k\tens s^m
 \end{align*}
 where in principle we sum over $k,m$. However, $s^is^k=0$ unless $k=i$ when $\eps_{ik}=0$, so the second term vanishes. The first term simplifies using $\extd (\sigma^i s^i)=-{\rm Vol}$ to the result stated. For the Ricci curvature, we use the standard lift $i(s\wedge t)={1\over 2}(s\tens t+t\tens s)$ and the inverse metric $(t,s)={1\over 2}=-(s,t)$, which translate to 
 \[ i({\rm Vol})={1\over 2\imath}(s^1\tens s^1+s^2\tens s^2),\quad (s^i,s^j)=-\imath\eps_{ij}\]
 respectively.  This then gives Ricci and $S$ as stated.  
 \endproof
 
 We see that there is a constant Ricci scalar and the Ricci tensor is proportional to the metric, much like a sphere. 
 
\subsection{Clifford action and spinor connection on $M_2(\C)$}\label{seclocM2}

Next we solve for the local tensorial level data for a geometrically realised Dirac operator for this metric. We take $J$ from Lemma~\ref{lemJ} as already analysed for a 2D calculus. Next, we analyse the content of the Clifford relations (\ref{cliffrel}) for the matrices $s^i\la=C^i$ with respect to a central spinor basis $\{e^\alpha\}$. Given the form of the wedge product and of the metric, the Clifford relations with $\varphi=\id$ and $\kappa=1$ then come down to 
 \begin{equation}\label{M2cliff} C^2 C^1=\imath \id, \quad C^1C^2+C^2C^1=0,\quad (C^1)^2=(C^2)^2.\end{equation}
The first two here have no joint solutions. However, we can make use of the weaker version with nontrivial automorphism $\varphi$.

\begin{lemma}\label{lemcliffmat}The full Clifford relations in the general form up to automorphism  hold if and only if $C^2C^1$ is invertible and the remainder  $C^1C^2+C^2C^1=0$, $(C^1)^2=(C^2)^2$ of (\ref{M2cliff}) hold. 
\end{lemma}
\proof  Given the relations of $\Omega^2$ in our case, and the off-diagonal form of the metric, the more general (\ref{cliffrel}) comes down in matrix terms in our basis to 
\[  C^2 C^1\varphi=\kappa g^{12}\id=-\imath\kappa,\quad C^1 C^2\varphi=\kappa g^{21}=\imath\kappa\id,\quad ((C^1)^2-(C^2)^2)\varphi=\kappa(g^{11}-g^{22})=0\]
 where $g^{ij}=(s^i,s^j)$ is the inverse matrix to the coefficient matrix $g_{ij}$ of the quantum metric. Without loss of generality, we can set $\kappa=1$ then we see that $\varphi=\imath(C^1C^2)^{-1}$ is uniquely determined and the only conditions then are that the latter parts of  (\ref{M2cliff}) hold. \endproof
 
 In the first instance, we will look more broadly by imposing only the last of (\ref{M2cliff}). Then one finds that there are principally 4 and 5 dimensional moduli spaces
 \[ {\rm I:}\quad  C^2=\pm C^1,\quad {\rm II:}\quad {\rm Tr}(C^1)={\rm Tr}(C^2)=0,\quad \det(C^1)=\det(C^2)\]
 We focus on the more interesting type II case, which we parameterise generically as
 \[ C^1= 
\begin{pmatrix}
c_1 & c_2 \\
 \frac{c_0-c_1^2}{c_2} & -c_1
\end{pmatrix}
,\quad C^2=\begin{pmatrix}
c_3 & c_4 \\
 \frac{c_0-c_3^2}{c_4} & -c_3
\end{pmatrix};\quad c_i\in \C, \]
where $-c_0$ is the common value of the determinant. There are also solutions with 0 in the top right corner but we focus on the above principal part of the moduli space for our detailed calculations. Type II includes choosing the $C^i$ as any two distinct  Pauli matrices. Next we solve the covariance of $\la$ conditions.

\begin{theorem}\label{thm_alt} For the $C^i$ of type (II) as above and any $J$ obeying (\ref{JJ}), we can solve all the remaining local tensorial level conditions (\ref{SJ}) -- (\ref{CJ}) and (\ref{dnablacliff}) for a fully geometrically realised spectral triple of inner type if and only if $\mu=0$ in the 1-parameter QLC. In this case, solutions have the form
\[  \sigma_\CS{}^i{}_j=\zeta^i{}_j [C^1,C^2], \quad \overline{C^i}=\eps' \zeta^i{}_j J [C^1,C^2]C^j J^{-1}\]
for any invertible matrix $\zeta\in M_2(\C)$  obeying
\[ \overline{\zeta}\zeta=\id,\quad \det(\zeta)= - {1\over \det([C^1,C^2])^2}.\]
Here $[C^1,C^2]$ is the matrix commutator. Moreover, a canonical $\gamma$ obeying its local conditions (\ref{gamsq})-(\ref{Sgam}) has the form 
  \[ \gamma={[C^1,C^2]\over \sqrt{-\det([C^1,C^2])}},\quad \overline{\gamma}=\eps'' J \gamma J^{-1}.\]
  These solutions are subject to solving the stated reality conditions for $C^i$ and $\gamma$ respectively for some $\eps',\eps''$. 
\end{theorem}
 \proof   For the $C^i$ of type (II) as above, we first solve  (\ref{dnablacliff}) for $\sigma_\CS$ using the equations for the inner case in Lemma~\ref{leminn} (the $S_i$ are given by (\ref{SthetaA}) with $A_i=0$) and  find the generic solution which we recognise as stated but provided $\mu=0$. We then insert this into the reality condition (\ref{CJ})  and solve the resulting condition by setting $\overline{C^i}$ as stated.  Finally, we put the form of $\sigma_\CS$ into (\ref{SJ}) in the inner form in Lemma~\ref{leminn} to get 
 \[  \delta_{jk}J=\overline{\zeta^k{}_j}\overline{[C^1,C^2]}J \zeta^i{}_j [C^1,C^2]. \]
 We solve this as two equations, one for $\overline{\zeta^i{}_j}$ as stated, which specifies the real and imaginary parts of $\zeta^i{}_j$, and the matrix equation 
 \[ J=\overline{[C^1,C^2]}J [C^1,C^2].\]
  We then solve the latter for $J$ of type (2)  in Lemma~\ref{lemJ} to obtain just the constraint shown for $\det(\zeta)$ and no constraint on the parameter $z$ of $J$. We then check that in fact any $J$ obeys this provided $\det(\zeta)$ is as shown. 
  
  For  $\gamma$, we put $\sigma_\CS$ as found into (\ref{gamsq}),(\ref{Cgam}) to find generically a unique solution, which we identify as stated. We take   
 (\ref{Jgam}) as a specification of $\overline{\gamma}$ as stated. We automatically obey (\ref{Sgam}) in the form (\ref{Sgaminn})  since $\sigma_\CS{}^i{}_j\propto [C^1,C^2]$ and this matrix commutes with $\gamma$ since the $C^i$ anticommute. \endproof
 
 This gives a canonical solution for inner geometrically realised spectral triples at the local tensorial level with 10 parameters: 5 for $C^i$, 3 for $\zeta$ and 2 for $J$ e.g. as in Lemma~\ref{lemJ}, modulo an overall unitary equivalence and for a choice of $\eps,\eps'$. These are subject to the reality constraints shown and we still need to solve for these. In the case of $\zeta$, the `reality condition' is the same as that of an antilinear isometry with respect to our basis and was already solved (with $\eps=1$) in Lemma~\ref{lemJ} when we solved for $J$, i.e. 
\[{\rm (1):}\quad  \zeta=e^{-{\imath\zeta_0}}\begin{pmatrix}
 \zeta_1 e^{\imath \zeta_3} & \zeta_2  \\
 \frac{1-\zeta_1^2}{\zeta_2} &- \zeta_1 e^{-\imath \zeta_3} 
\end{pmatrix},\quad {\rm (2):}\quad \zeta=e^{-{\imath\zeta_0}}\begin{pmatrix}
 e^{\imath \zeta_2} & 0  \\
 \zeta_1  &- e^{-\imath \zeta_2} 
\end{pmatrix};\quad \zeta_i\in \R\]
in polar coordinate form with real parameters. We have explicitly included a phase factor at the front whereby $\det(\zeta)=-e^{-2\imath \zeta_0}$. For each $\zeta$, we still have to solve for $C^i$ with 
\[\det([C^1,C^2])=4 c_0^2-\frac{\left((c_2 c_3-c_1 c_4)^2-c_0 \left(c_2^2+c_4^2\right)\right)^2}{c_2^2 c_4^2}=e^{\imath\zeta_0}\]
and obeying the reality condition in the theorem. The system of equations for the latter is rather nonlinear and we do not give the most general solutions. Rather, we  focus on solutions where the entries of $J$ and $C^i$ are real and some obvious choices of $\zeta$. 

\begin{example}\label{ex_alt} We set $\zeta=\sigma^3$ and $J$ to be real of type (1), $C^i$ with real entries and the simplifying assumption that $c_2=c_4=1$. There are then no solutions with $\eps=-1$, while for $\eps=1$ we have:

(a)  A 1-parameter family of  solutions with $\eps'=1$, namely
\[ J=\begin{pmatrix}x & 2 \\
 \frac{1-x^2}{2} & -x\end{pmatrix},\quad C^1=\begin{pmatrix}
 \frac{x+1}{2} & 1 \\
 -({x+1\over 2})^2 & -\frac{x+1}{2} \end{pmatrix},\quad C^2=\begin{pmatrix}
 \frac{x-1}{2} & 1 \\
 - ({x-1\over 2})^2 & \frac{1-x}{2}\end{pmatrix}\]
 with a single real parameter $x\in \R$. Here $\det(C^i)=0$, while $\det([C^1,C^2])=-1$ as required, and the relations
 \[ C^1C^2+C^2C^1=-\id.\]
 Moreover, the canonical form of $\gamma$ works and is given by
 \[ \gamma=[C^1,C^2]=J,\quad \eps''=1.\]
 We have similar solutions with the $C^i$ swapped and/or with the off-diagonals of $J$ negated and $C^i$ negated. 
 
 (b) A 1-parameter family with $\eps'=-1$, namely the same $J$ as before and 
 \[  C^1=\begin{pmatrix}
 \frac{x+1}{2} & 1 \\
 \frac{1-x (x+2)}{4} & -\frac{x+1}{2}\end{pmatrix},\quad  C^2=\begin{pmatrix}
 \frac{x-1}{2} & 1 \\
 \frac{1-x(x-2)}{4} & \frac{1-x}{2} \end{pmatrix}\]with $\det(C^i)=-{1\over 2}$. These have 
\[ C^1C^2+C^2C^1=0\]
so that the full Clifford relations hold in a generalised form with automorphism. 
Moreover, the canonical form of $\gamma$ works and is given by
 \[ \gamma=-\imath [C^1,C^2]=-\imath\begin{pmatrix}x & 2 \\
 -\frac{1+x^2}{2} & -x\end{pmatrix},\quad \eps''=1.\]
There are similarly 3 variants as in case (a). Also note that  $J$ in the two cases can never be symmetric. 
\end{example}

\begin{remark} The example of a geometrically realised spectral triple on $M_2(\C)$ in \cite{BegMa:spe,BegMa} has the same calculus and quantum metric as used above. Working with our self-adjoint basis $s^i$ and taking $e^1=e,e^2=f$ in our conventions, one has $\eps=\eps'=-1$ and
\[  J=\begin{pmatrix} 0 & -1 \\ 1 & 0\end{pmatrix}=-\imath\sigma^2,\quad C^1=-\imath\sigma^1,\quad C^2=\imath\sigma^2, \quad [C^1,C^2]=\imath\sigma^3,\quad \sigma_\CS{}^i{}_j=\delta^i{}_j\id.\]
Here $C^1C^2+C^2C^1=0$ and $(C^1)^2=(C^2)^2$ (both being $-\id$), so that the full Clifford relations hold with a suitable automorphism $\varphi\propto \sigma^3$. But this example fundamentally did not have covariance of $\la$ i.e. (\ref{SJ}) fails as $S_i=0$ while the $\mu=0$ QLC does not have zero Christoffel symbols. 
By contrast, we see that we can have solutions with $\eps=1, \eps'=-1$ with covariance of $\la$ as well as the full Clifford relations generalised to allow an automorphism $\varphi$. The difference is that $\zeta=\sigma_3$ and $\sigma_\CS$ also has a nontrivial operation $[C^1,C^2]$ as per Theorem~\ref{thm_alt} rather than the identity on the spinor indices . 
\end{remark}

The other obvious real choice of $\zeta$ also has solutions as follows.

\begin{example}\label{ex_altp} We set $\zeta=\sigma^1$ and $J$ to be real of type (1)   and $C^i$ real with  $c_2=c_4=1$. There are then no solutions  with $\eps=-1$, while for $\eps=1$  we find have: 

(a) A  2-parameter  family, for either sign of $\eps'$, namely
\[ J=
\left(
\begin{array}{cc}
 \frac{\eps'+x^2-y^2}{\sqrt{2}} & {\scriptstyle \sqrt{2} (x-y)} \\
 \frac{1-\frac{1}{2} \left(\eps'+x^2-y^2\right)^2}{\sqrt{2} (x-y)} & -\frac{\eps'+x^2-y^2}{\sqrt{2}} \\
\end{array}
\right),\]
 \[ C^1=\left(
\begin{array}{cc}
 x & 1 \\
\frac{(x-y)^4+1}{4 (x-y)^2}-x^2\quad &  -x \\
\end{array}
\right),\quad C^2=\left(
\begin{array}{cc}
 y & 1 \\
 \frac{(x-y)^4+1}{4 (x-y)^2}-y^2\quad & -y \\
\end{array}
\right)\]
with
\[ C^1C^2+C^2C^1=\frac{1}{2} \left(\frac{1}{(x-y)^2}-(x-y)^2\right)\id,\quad \det(C^2C^1)=\frac{\left((x-y)^4+1\right)^2}{16 (x-y)^4}. \]
Hence, if we set $x-y=\pm1$ then $C^1C^2+C^2C^1=0$ and $C^2C^1$ is invertible, i.e. the full Clifford relations hold with automorphism. Moreover, the canonical form of $\gamma$ works and is given by
 \[ \gamma=-\imath [C^1,C^2]=-\imath\begin{pmatrix}
 x^2-y^2 & 2 (x-y) \\
 -\frac{\left(x^2-y^2\right)^2+1}{2 (x-y)} & y^2-x^2  \end{pmatrix},\quad \eps''=1.\]
There is a similar solution with the opposite sign of $\sqrt{2}$.  Note also that one can choose $x,y$ such that $J$ is symmetric but not at the points where $x-y=\pm 1$. 

(b) There is a 1-parameter family  for  $\eps'=-1$ only, namely 
\[ J=\left(
\begin{array}{cc}
 0 & {\scriptstyle \sqrt{4 x \left(x-\sqrt{x^2-1}\right)-2}} \\
 \frac{1}{\sqrt{4 x \left(x-\sqrt{x^2-1}\right)-2}} & 0 \\
\end{array}
\right),\quad C^1=\left(
\begin{array}{cc}
 x & 1 \\
 -\frac{1}{2} & -x \\
\end{array}
\right),\quad C^2= \left(
\begin{array}{cc}
 \sqrt{x^2-1} & 1 \\
 \frac{1}{2} & -\sqrt{x^2-1} \\
\end{array}
\right) \]
with
\[ C^1C^2+C^2C^1=2 x \sqrt{x^2-1}\id,\quad \det(C^2C^1)={1\over 4}+x^2(x^2-1).\]
Here $|x|\ge 1$, so the $C^i$ anticommute only at $x=\pm 1$. Moreover, the canonical form of $\gamma$ works and is given by
 \[ \gamma=-\imath [C^1,C^2]=-\imath\begin{pmatrix}
 1 & 2 (x- \sqrt{x^2-1}) \\
 -\sqrt{x^2-1}-x & -1 \end{pmatrix},\quad \eps''=1.\]
There are 3 similar solutions with different signs for the two types of square roots independently.  Note also that $J$ is symmetric at $x={3\over 2\sqrt{2}}$. 

(c) A 1-parameter family of solutions, for either sign of $\eps'$, namely
\[ J=\left(
\begin{array}{cc}
 \frac{1}{\sqrt{2}} &{\scriptstyle \eps' 2 \sqrt{2} x }\\
 \frac{\eps'}{4 \sqrt{2} x} & -\frac{1}{\sqrt{2}} \\
\end{array}
\right),\quad C^1=\left(
\begin{array}{cc}
 x & 1 \\
 \frac{1}{16 x^2} & -x \\
\end{array}
\right),\quad C^2=\left(
\begin{array}{cc}
 -x & 1 \\
 \frac{1}{16 x^2} & x \\
\end{array}
\right)\]
\[ C^1C^2+C^2C^1=({1\over 8 x^2} - 2 x^2)\id,\quad \det(C^2C^1)={1\over 8}+ x^4+ {1\over 256 x^4}.\]
We see that $C^1C^2+C^2C^1=0$ at $x=\pm {1\over 2}$ and $C^2C^1$ is always invertible, so in  this case we have the full Clifford relations with automorphism. Moreover, the canonical form of $\gamma$ works and is given by
 \[ \gamma=-\imath [C^1,C^2]=-\imath\begin{pmatrix}
 0 & 4 x \\
 -\frac{1}{4 x} & 0 \end{pmatrix},\quad \eps''=1.\]
There is a similar solution with the opposite sign of $\sqrt{2}$. Note also that $J$ is symmetric at $x=\pm{1\over 4}$. 
\end{example}

We can equally let $J$ be of type (2), again restricting to real entries at least for the $C^i$. The algebra is easier and we do not need to impose further conditions on $c_2,c_4$. 

\begin{example}\label{exJ2} We set $\zeta=\sigma^3$ and $J$ to be type (2)  with $\eps=1$ and real entries for $C^i$.  Then for either sign of $\eps'$, we have two 2-parameter solutions: 
\[ J=\left(
\begin{array}{cc}
 1 & 0 \\
 -{2 x\pm1\over y} & -1 \\
\end{array}
\right),\quad C^1=\left(
\begin{array}{cc}
 x & y \\
{ \frac{1}{2}-x^2\over y} & -x \\
\end{array}
\right),\quad C^2=\eps'\left(
\begin{array}{cc}
x\pm1 &   y \\
{\frac{1}{2}-(x\pm1)^2\over y}\  & -(x\pm1 )\\
\end{array}
\right),\]
\[ C^1C^2+C^2C^1=0,\quad \gamma=-\imath [C^1,C^2]=-\imath\eps'\left(
\begin{array}{cc}
 -1\mp 2 x & \mp2 y \\
 \pm{2 x (x\pm 1)+1\over y} &1 \pm 2 x \\
\end{array}
\right)\]
with $\eps''=1$. Here $C^1C^2$ is invertible and hence this obeys the full Clifford relations in the general form with automorphism. Note that $J$ is symmetric at $x=\mp {1\over 2}$ and that the $y$ parameter can be absorbed in a conjugation and hence we could set $y=c_2=c_4=1$ as we did before. There are also some 1-parameter solutions with only this $y$ parameter, i.e. discrete if we ignore this freedom. The matrix entries of $J$ here are also real but allowing the full complex case does not lead to fundamentally new solutions. 
 \end{example}

We find no solutions for $\zeta=\sigma^1$ under similar assumptions. 
It should be clear that there are potentially many other solutions, because we do not need to asssume $c_2=c_4=1$ in the $J$ type (1) examples and we do not need to restrict to real entries for $C^i$ or $J$. These restrictions are not too impactful given the freedom (\ref{unitaryequiv}) to conjugate solutions. Also note that real entries  holds classically in two dimensions if one uses $\sigma^1,\sigma^3$ for the Clifford algebra rather than the more standard $\sigma^1,\sigma^2$.  

\subsection{Natural spectral triple on $M_2(\C)$}\label{secspecM2}

Looking at the solutions in Section~\ref{seclocM2} for a geometrically realised spectral triple at the local tensorial level (i.e. before considering the Hilbert space structure), all solutions under our assumptions have $\eps=1$. Moreover, the only ones for which the $C^1,C^2$ anticommute (and the full Clifford relations hold with automorphism) and $J$ is symmetric  (so that $\CJ$ is an antilinear isometry with respect to the usual trace inner product given by $\mu^{\alpha\beta}=\delta^{\alpha\beta}$ for the measure) are $x=\pm {1\over 2}$ in Example~\ref{exJ2}.  This is therefore the natural solution for a full geometric realisation under our simplifying assumptions, mainly of real entries in the $C^i$. The value of $y$ is not significant as it can be changed by conjugation. 

To be concrete, henceforth we focus in particular on $x=-{1\over 2}$ and  $y={1\over 2}$  in Example~\ref{exJ2} so that
\begin{equation}\label{canJtype2} J=\sigma^3=\zeta,\quad C^1={\sigma^1-\sigma^3\over 2},\quad C^2=\eps'{\sigma^1+\sigma^3\over 2},\quad \gamma=-\eps'\sigma^2,\quad \eps=1=\eps''\end{equation}
in terms of Pauli matrices.  Then $\CJ$ is an antilinear isometry. We also need $\gamma$ to be hermitian. As an operator, $\gamma(e^\alpha)=\imath\eps'\eps_{\alpha\beta}e^\beta$ (with sum over repeated indices), where $\eps_{\alpha\beta}$ is the antisymmetric tensor with $\eps_{12}=1$. Hence on the spinor components of $\psi=\psi_\alpha e^\alpha$, we have likewise
\begin{equation}\label{gammasol}(\gamma\psi)_\beta=-\eps'\psi_\alpha \sigma^2{}^\alpha{}_\beta=\imath\eps' \psi_\alpha \eps_{\alpha\beta}.\end{equation} Then  
\[ \langle\gamma\phi,\psi\rangle={\rm Tr}(\gamma\phi)_\beta^\dagger\psi_\beta={\rm Tr} (\imath\eps'\eps_{\alpha\beta}\phi_\alpha)^\dagger\psi_\beta={\rm Tr}\phi_\alpha^\dagger(-\imath\eps' \eps_{\alpha\beta}\psi_\beta)={\rm Tr}\phi_\alpha^\dagger(\gamma\psi)_\alpha=\langle\phi,\gamma\psi\rangle\]
as required.

\begin{proposition}\label{propspecM2} The solution (\ref{canJtype2}) gives a full geometrically realised spectral triple which is  $n=0$ even if $\eps'=-1$ with
\[ (\dirac\psi)_\beta=-{\imath\over 4}\left( [-\sigma^1+\sigma^2,\psi_\alpha]\sigma^1{}^\alpha{}_\beta+\{\sigma^1+\sigma^2,\psi_\alpha\} \sigma^3{}^\alpha{}_\beta    \right)\]
 and  $n=1$ ignoring $\gamma$ if $\eps'=1$ with
\[ (\dirac\psi)_\beta={\imath\over 4}\left( \{\sigma^1+\sigma^2,\psi_\alpha\}\sigma^1{}^\alpha{}_\beta+[-\sigma^1+\sigma^2,\psi_\alpha] \sigma^3{}^\alpha{}_\beta    \right).\]
\end{proposition}
\proof We use (\ref{dirac}) and  the form of $\del_i$ in (\ref{M2deli}) and of  $S^\alpha{}_{i\beta}$ in (\ref{SthetaA}) gives the Dirac operator
 \begin{equation}\label{diracinner} (\dirac \psi)_\gamma={\imath \over 2}(\sigma^i\psi_\alpha-\psi_\beta\sigma^j \sigma_\CS{}^{\beta j}{}_{i\alpha})C^{i\alpha}{}_\gamma.\end{equation}
For the chosen solution and its form of $\sigma_\CS{}^i{}_j$, and noting that its $C^i$ obey
\[ [C^1,C^2]C^1=-C^2,\quad [C^1,C^2]C^2= C^1,\]
we have
\begin{align}\label{diracsol}  (\dirac \psi)_\beta&={\imath\over 4}\left( (\sigma^1\psi_\alpha+\psi_\alpha \sigma^2)(\sigma^1-\sigma^3)^\alpha{}_\beta+\eps'(\sigma^2\psi_\alpha+\psi_\alpha\sigma^1)(\sigma^1+\sigma^3)^\alpha{}_\beta \right)\end{align} which, collecting terms, we can write as stated.  Note that the $\sigma^1,\sigma^2$ coming from $\theta$ are simply elements of the algebra while the $\sigma^1,\sigma^3$ coming from the $C^1,C^2$ etc. act as shown explicitly on spinor indices. Hence
\begin{align*} 
\< \phi,\dirac\psi\>&={\imath\over 4}{\rm Tr} \phi_\alpha^\dagger\left( (\sigma^1\psi_\beta+\psi_\beta \sigma^2)(\sigma^1-\sigma^3)^\beta{}_\alpha+\eps'(\sigma^2\psi_\beta+\psi_\beta\sigma^1)(\sigma^1+\sigma^3)^\beta{}_\alpha \right)\\
\<\dirac\phi,\psi\>&=-{\imath\over 4}{\rm Tr} \left( (\sigma^1\phi_\alpha+\phi_\alpha \sigma^2)(\sigma^1-\sigma^3)^\alpha{}_\beta+\eps'(\sigma^2\phi_\alpha+\phi_\alpha\sigma^1)(\sigma^1+\sigma^3)^\alpha{}_\beta \right)^\dagger\psi_\beta\\
&=-{\imath\over 4}{\rm Tr} \left( (\phi_\alpha^\dagger \sigma^1+\sigma^2\phi_\alpha^\dagger)\psi_\beta(\sigma^1-\sigma^3)^\beta{}_\alpha+\eps'(\phi_\alpha^\dagger \sigma^2+\sigma^1\phi_\alpha^\dagger)\psi_\beta(\sigma^1+\sigma^3)^\beta{}_\alpha \right)\\
&=-\<\phi,\dirac\psi\>
\end{align*}
so that $\dirac$ is antihermitian. In the second computation, we used that $\sigma^1\pm\sigma^3$ are real to move them out from under $(\ )^\dagger$ and that they are symmetric to swap $\alpha,\beta$. We also used that $\sigma^1,\sigma^2\in M_2(\C)$ are hermitian as elements of the algebra. We then compare the terms using cyclicity of the trace to  obtain $-\<\phi,\dirac\psi\>$. 

Finally, remembering that the Dirac operator in a Connes spectral triple is hermitian and hence corresponds to $D=\imath\dirac$, the interpretation of $\eps'$ is reversed. Hence our above construction is $n=0$ in the period 8 classification if $\eps'=-1$, as an even spectral triple, and  $n=1$ ignoring $\gamma$, as an odd spectral triple if $\eps'=1$. \endproof

One can check directly that the stated $D=\imath\dirac$, $\gamma$ in (\ref{gammasol}) and 
\begin{equation}\label{Jsol} (\CJ\psi)_\beta=\psi_\alpha^\dagger\sigma^3{}^\alpha{}_\beta\end{equation}
 indeed obeys all the required  axioms for a Connes spectral triple. Here 
\[ (\CJ b\CJ^{-1}\psi)_\beta=\psi_\beta b^\dagger,\quad ([D,a]\psi)_\beta={1\over 4}\left([\sigma^1-\sigma^2,a]\psi_\alpha \sigma^1{}^\alpha{}_\beta+[\sigma^1+\sigma^2,a]\psi_\alpha\sigma^3{}^\alpha{}_\beta\right)\]
in the even case, and a similar formula with $\sigma^1,\sigma^3$ swapped and ignoring $\gamma$ in the odd case. Moreover, underlying $D$ is  the bimodule connection $\nabla_\CS$ which for our basis has Christoffel symbols
\begin{equation}\label{Sisol} S_i^\alpha{}_\beta={\imath\over 2}\left(\sigma^i\delta^\alpha{}_\beta-\sigma^j\sigma^3{}^j{}_i[C^1,C^2]^\alpha{}_\beta\right)\end{equation}
which in our case is 
\begin{equation}\label{S1S2mat} S{}_1={\imath\over 2}\sigma^1\begin{pmatrix}1 & \eps'\\ -\eps' &1\end{pmatrix},\  S{}_2={\imath\over 2}\sigma^2\begin{pmatrix}1 & -\eps'\\ \eps' &1\end{pmatrix}.\end{equation}
We have seen that this arises naturally up to conjugation under some modest assumptions and the requirement of full quantum-geometric realisability. 

\begin{corollary}\label{corM2dirac} The above solution (\ref{canJtype2}) for the geometric data is unitarily equivalent to
 \[ C^1={1\over\sqrt{2}}\sigma^1,\quad C^2 ={\eps'\over \sqrt{2}}\sigma^3,\quad J={1\over \sqrt{2}}(\sigma^3-\sigma^1)\]
 with $\gamma,\sigma^i_S{}_j,S_i$ unchanged and gives $n=0,1$ spectral triples with
 \[ (\dirac \psi)_\beta={\imath\over 2 \sqrt{2}}\left( (\sigma^1\psi_\alpha+\psi_\alpha \sigma^2)\sigma^1{}^\alpha{}_\beta+\eps'(\psi_\alpha\sigma^1+\sigma^2\psi_\alpha)\sigma^3{}^\alpha{}_\beta \right).\]
 \end{corollary}
 \proof Here we apply (\ref{unitaryequiv}), i.e. a change of the $\{e^\alpha\}$ basis, via the orthogonal real matrix
\[ u=\begin{pmatrix}
 \frac{\sqrt{2+\sqrt{2}}}{2} & \frac{\sqrt{2-\sqrt{2}}}{2} \\
 -\frac{\sqrt{2-\sqrt{2}}}{2}  & \frac{\sqrt{2+\sqrt{2}}}{2} 
 \end{pmatrix}\]
to obtain new $C^i,J$ as stated and $u\dirac u^{-1}$ for $\dirac$ in (\ref{diracsol}) to obtain the new Dirac operator stated. Here $[C^1,C^2]=-\imath\eps'\sigma^2$ is not changed by the conjugation so $\gamma,\sigma^i_S{}_j,S_i$ are not affected. As a check, one obtains the same new $\dirac$ from (\ref{dirac}) using $\del_i$ from (\ref{M2deli}),  the new $C^i,J$ and the connection $\nabla_\CS$ with $S_i$ given by (\ref{Sisol}) and $\gamma$ by the formula in Theorem~\ref{thm_alt}. \endproof

As a further check on all of our calculations, we verify that the Lichnerowicz formula (\ref{diracsqmat}) holds by computing the relevant quantities directly for our above quantum geometry.

\begin{proposition}\label{proplichM2} The Dirac operator $\dirac$ in Corollary~\ref{corM2dirac} and the underlying quantum geometry has
\[ (\dirac^2\psi)_\beta=-{1\over 2}\psi_\beta  - {1\over 4}(\sigma^1\psi_\beta\sigma^2+\sigma^2\psi_\beta\sigma^1)+ {\eps'\over 4}\{\sigma^3,\psi_\alpha\}\sigma^{2\alpha}{}_\beta,\quad  \varphi^\alpha{}_\beta=-2\eps'\sigma^{2\alpha}{}_\beta, \]
\[ (\square_\CS\psi)_\beta= -{1\over 2}\{\sigma^3,\psi_\beta\}+{\eps'\over 2}(\sigma^2\psi_\alpha\sigma^1+\sigma^1\psi_\alpha\sigma^2),\quad R_\CS{}^\alpha{}_\beta=\eps'\sigma^{2\alpha}{}_\beta\]
and $\kappa=1$, in agreement with (\ref{diracsqmat}).
\end{proposition}
\proof We let
\[ \tilde\psi_\alpha:= (\sigma^1\psi_\gamma+\psi_\gamma\sigma^2)\sigma^{1\gamma}{}_\alpha+\eps'(\psi_\gamma\sigma^1+\sigma^2\psi_\gamma)\sigma^{3\gamma}{}_\alpha\]
as  in $(\dirac\psi)_\beta$ and put this into
\begin{align*} -8(\dirac^2\psi)_\beta&=  (\sigma^1\tilde\psi_\alpha+\tilde\psi_\alpha\sigma^2)\sigma^{1\alpha}{}_\beta+\eps'(
\tilde\psi_\alpha\sigma^1+\sigma^2\tilde\psi_\alpha)\sigma^{3\alpha}{}_\beta,
\end{align*}
which expands out to give $\dirac^2$ as stated on collecting terms with and without $\eps'$. We used the Pauli algebra relations $\sigma^1\sigma^2=\imath\sigma^3=-\sigma^2\sigma^1$ and $(\sigma^1)^2=(\sigma^2)^2=1$ in $A$. At the same time we used $(\sigma^1\sigma^3)^\alpha{}_\beta=-\imath\sigma^{2\alpha}{}_\beta=-(\sigma^3\sigma^1)^\alpha{}_\beta$ and again that $\sigma^1,\sigma^3$ square to 1, now in the matrix algebra with spinor indices. Similarly, from the proof of Lemma~\ref{lemcliffmat} and $C^i$ from Corollary~\ref{corM2dirac}, we have $\kappa=1$ and $\varphi=\imath (C^1C^2)^{-1}=\imath 2\eps'(\sigma^1\sigma^3)^{-1}=-2\eps'\sigma^2$ as stated.  

For $R_\CS{}^\alpha{}_\beta$, we use Lemma~\ref{lemRS} and   
\[ (C^{j\beta}{}_\gamma C^{i\gamma}{}_\delta \varphi^\delta{}_\eta -\kappa g^{ij} \delta^\beta{}_\eta)={1\over 2}\varphi^\beta{}_\eta\delta^{ij}.\]
Here, for $i=1,j=2$, this expression is $C^2C^1\varphi - (-\imath)\id=-C^1C^2\varphi+\imath\id=0$ from the definition of $\varphi$. Similarly for the other cases. Also, our QLC Christoffel symbols and braiding from Section~\ref{secM2qrg} for the $\mu=0$ QLC $\nabla s^i=2\theta\tens s^i$  are 
\[ -{1\over 2}\Gamma^k{}_{ij}=\imath \sigma^i \delta^k{}_j,\quad \sigma^{ij}{}_{kl}=-\delta^i{}_l\delta^j{}_k\]
and we recall that $\del_i={\imath\over 2}[\sigma^i,\ ]$ for our calculus. Then  
\begin{align*}R_\CS{}^\alpha{}_\eta&=-\eps'({1\over 2}[\sigma^i,S^\alpha{}_{i\beta}]+S^\alpha{}_{i\beta}\imath\sigma^i-S^\alpha{}_{i\gamma} S^\gamma{}_{i\beta})\sigma^{2\beta}{}_\eta\\
&=-\eps'({1\over 2}\{\sigma^i,S^\alpha{}_{i\beta}\}-S^\alpha{}_{i\gamma} S^\gamma{}_{i\beta})\sigma^{2\beta}{}_\eta\\
&=-{\eps'\over 4}\big(-2(\id+\imath\eps'\sigma^2)- 2(\id-\imath\eps'\sigma^2)  +(\id+\imath\eps'\sigma^2)^2-(\id-\imath\eps'\sigma^2)^2\big)^\alpha{}_\beta \sigma^{2\beta}{}_\eta \end{align*}
on substituting $S^\alpha{}_{i\beta}={\imath\over 2}\sigma^i(\id\pm \eps'\sigma^2)^\alpha{}_\beta$ for $i=1,2$ respectively from (\ref{S1S2mat}).  This then simplifies to the answer stated after cancellations. 

Next, for the Laplacian on spinors, we use Proposition~\ref{propLapS}, first computing 
\begin{align*} L^\alpha{}_\beta&= -\imath\del_1 S^\alpha{}_{2\beta}+\imath \del_2S^\alpha{}_{1\beta}- S^\alpha{}_{1\beta}\sigma^2+ S^\alpha{}_{2\beta}\sigma^1- \imath S^\alpha{}_{1\gamma} S^\gamma{}_{2\beta}+\imath S^\alpha{}_{2\gamma} S^\gamma{}_{1\beta}\\
&={1\over 2}\{\sigma^1,S^\alpha{}_{2\beta}\} - {1\over 2}\{\sigma^2, S^\alpha{}_{1\beta}\} - \imath S^\alpha{}_{1\gamma} S^\gamma{}_{2\beta}+\imath S^\alpha{}_{2\gamma} S^\gamma{}_{1\beta}=-\sigma^3\delta^\alpha{}_\beta
\end{align*}
where we use $-{1\over 2}g^{kl}\Gamma^j{}_{kl}= -\sigma^2, \sigma^1$  for $j=1,2$ respectively and $g^{kl}\sigma^{ij}{}_{kl}=-g^{ji}=g^{ij}$  to obtain the first expression. The form of $\del_i$ gives the second expression and then substituting  the values of $S_i$ as above and simplifying using the Pauli algebra relations gives the final result. Similarly, if $a\in A$, the scalar Laplacian from Proposition~\ref{propLapS} is 
\[ \square a=\imath\del_2\del_1 a-\imath\del_1\del_2 a+ (\del_1 a)(-\sigma^2)+(\del_2 a)\sigma^1=-{1\over 2}\{\sigma^3,a\}+{\imath \over 2}(\sigma^2 a\sigma^1-\sigma^1 a \sigma^2).\]
(This has $\square \sigma^3=-2$ and is otherwise zero on the basis $1,\sigma^i$.) Putting all this together, we have 
\[(\square_\CS\psi)_\beta=-{1\over 2}\{\sigma^3,\psi_\beta\}+{\imath \over 2}(\sigma^2 \psi_\beta \sigma^1-\sigma^1 \psi_\beta  \sigma^2) -\psi_\beta\sigma^3- 2\imath(\del_1\psi_\alpha)S^\alpha{}_{2\beta}+2\imath (\del_2\psi_\alpha)S^\alpha{}_{1\beta}\]
into which we put the form of $\del_i$ and $S_j$ and simplify, to give the result stated. From the final expressions stated, it is then easy to see that  (\ref{diracsqmat}) holds.  
\endproof

We see that even though the QLC on $M_2(\C)$ was flat, the spinor connection is not flat and moreover the curvature operator acts by multiplication by the constant value $-\sigma^3\in A$. We also see that the spinor Laplacian has a simple form on the components. 

 \section{Concluding remarks}\label{secrem}

In this paper, we set up the formalism for geometric realisation of spectral triples with central bases for $\Omega^1$ and for the spinor bundle $\CS$ as in  \cite{LirMa2} but now extended to allow the Christoffel symbols to be nonconstant and the basic 1-forms to not form a Grassmann algebra. We demonstrated how the formalism leads uniquely (up to unitary equivalence) to the expected $n=2$ spectral triple on the noncommutative torus, just as we did for the fuzzy sphere in \cite{LirMa2}, as well as $n=1$ for another choice of signs. Interestingly, we saw in Section~\ref{sectorext} that the spinor connection $\nabla_\CS$ and the connection $\nabla$ on $\Omega^1$ need not be trivial and can be used to encode a dynamically generated mass term or chiral term for the Dirac operator. It would be interesting to see how our results extend to a noncommutative $d$-torus, where it is known that there are $2^d$ spin structures matching their classical counterparts\cite{Ven}. Moreover, the noncommutative torus in a functional analysis context admits genuinely projective modules\cite{Con:dif}, and some of these could potentially lead to more sophisticated examples using the quantum geometric formalism of \cite{BegMa:spe}.

Looking beyond standard noncommutative torus, we note that the $q$-deformed noncommutative torus $\C_{q,\theta}[\T^2]$ is the same algebra as before but now with a $q$-deformed differential structure\cite[Exercise~1.5]{BegMa},
 \[( \extd u )u=q u\extd u,\quad (\extd v)v=qv\extd v,\quad (\extd v)u=e^{\imath\theta}u\extd v,\quad (\extd u)v=e^{-\imath\theta}v\extd u\]
 with $q$ real for a $*$-calculus.  One might ask if one can extend the analysis we have done to construct spectral triples in this case.  
The problem here is the lack of a central self-adjoint basis, which makes calculations harder. On the other hand, the $q$-deformed case is known to be a comodule algebra cocycle twist of the standard noncommutative torus\cite[Exercise~9.7]{BegMa} and hence can be constructed as a twist of the standard one above. Here, the quantum Riemannian geometry itself can be twisted using \cite{BegMa:twi}. 

Alternatively, we can set $q=e^{\imath\theta}$ and then $s^1=v^{-1}u^{-1}\extd u$ and $s^1=uv^{-1}\extd v$ are a central basis \cite[Exercise 1.5]{BegMa} but, because $q$ is not real, this case is not a $*$-calculus so we leave
the existing framework\cite{BegMa:spe} for geometric realisation of spectral triples. However, this case is nevertheless interesting and if 
$q=e^{\imath\theta}$ is a primitive $n$-th root of unity then one can define the reduced noncommutative torus as
\[ c_q[\T^2]= \C_{q,\theta}[\T^2]/ \< u^n-1, v^n-1\>,\]
which is compatible with the differential calculus \cite[Exercise1.5]{BegMa}, just not as a $*$-calculus. Just as $c_q[SL_2]$ are not usual Hopf $*$-algebras but `flip' ones, one could expect a theory for the reduced noncommutative torus that covers this non-standard case and possibly geometrically realises the spectral triples on fuzzy tori in \cite{BG}. This is another interesting direction for further work.  

Secondly, we demonstrated the formalism on $M_2(\C)$ with its standard 2D calculus and the two choices of quantum metric for which the QRG is known. The moduli of geometric Dirac operators at the local tensorial level turned out to be very rich and for the main quantum metric the Hilbert space conditions led naturally to a pair of spectral triples depending on $\eps'$. The other metric was relegated to the appendix since, as one might expect from the Lorentzian form of the metric, the most reasonable local Dirac operator was not antihermitian so that we do not have exactly a spectral triple. Indeed, one might expect a Lorentzian spectral triple \cite{Dev, PasSit}, although our example does not appear to immediately fit into one of these frameworks either. This is another topic for further work. 

Quantum Riemannian geometries have been used for baby quantum gravity models, see e.g. \cite{LirMa,ArgMa4} and the inclusion of a Dirac operator should allow the coupling of such models to fermionic matter, which remains to be investigated. There are also new tools such as quantum geodesics\cite{BegMa:cur}  recently computed on the noncommutative torus and $M_2(\C)$. Another reason to be interested in a full classification of geometrically realised spectral triples on finite-dimensional algebras such as $M_2(\C)$ and $c_q[\T^2]$ at roots of unity is the potential for applications to the structure of elementary particles by tensoring spacetime with such a finite quantum geometry\cite{Con0,ConMar}. In the current state of the art, while one can choose the algebra, its representation and $\dirac$ to reproduce the standard model and obtain  some predictions, having further (geometric realisation) constraints on the allowed Dirac operators would translate directly into further constraints on the particle physics. One can also come at the problem in a Kaluza-Klein manner without the Dirac operator\cite{ArgMa4,LiuMa2} and then extend this to spinors and spectral triples. This will be looked at elsewhere.

\appendix 
\section{Dirac operator for alternate QRG of $M_2(\C)$}\label{secM2alt}   

In this appendix we repeat the analysis as in Section~\ref{secM2} with the same $\Omega^1$ and self-adjoint basis $s^i$, but now for the alternate quantum metric  \cite[Exercise~8.5]{BegMa}, which is the only other one for which the QLCs have been analysed. Here $\dirac$ turns out not to be (anti)hermitian, unlike for the metric studied in the main text. 

\subsection{Alternate quantum Riemannian geometry on $M_2(\C)$.} \label{secM2qrgalt} Here the alternate metric (times 2) is $2(s\tens s+t\tens t)$, which in the self-adjoint basis comes out as
\begin{equation}\label{M2g_alt} g=-s^1\tens s^1 +  s^2\tens s^2,\end{equation}
with a `Lorentzian' metric.  The most natural 1-parameter QLC here (not the most general, but sufficiently nontrivial for our purposes) is
from \cite[Exercise~8.3]{BegMa},
\begin{align*}\label{nablast}\nabla s=2E_{21} t\tens s&-\rho E_{12}(s\tens t-t\tens s)+\rho E_{21}(t\tens t-s\tens s),\\
\nabla t=2E_{12}s\tens t&-\rho E_{21}(s\tens t-t\tens s)+\rho E_{12}(t\tens t-s\tens s),\end{align*}
where $\rho$ is an {\em imaginary} parameter.  This is again of an inner form determined by 
\[ \sigma(s\tens s)=s\tens s+\rho(s\tens t-t\tens s),\quad \sigma(t\tens t)=t\tens t+\rho(s\tens t-t\tens s),\]
\[\sigma(s\tens t)=-t\tens s+\rho(s\tens s-t\tens t),\quad \sigma(t\tens s)=-s\tens t+\rho(s\tens s-t\tens t).\]
The braiding converts to
\begin{align*} \sigma(s^1\tens s^1)&=-s^2\tens s^2- 2\imath\rho s^2\tens s^1,\quad \sigma(s^1\tens s^2)=s^1\tens s^2,\\
\sigma(s^2\tens s^1)&=s^2\tens s^1,\quad \sigma(s^2\tens s^2)=-s^1\tens s^1 - 2\imath\rho s^1\tens s^2,\end{align*}
which we can write as 
\begin{equation}\label{M2sigalt} \sigma(s^i\tens s^j)=\begin{cases}s^i\tens s^j & {\rm if\ }i\ne j\\ -s^{\bar i}\tens s^{\bar i}-2\imath\rho s^{\bar i}\tens s^i & {\rm if\ }i=j,\end{cases}\end{equation}
where $\bar 1=2$ and $\bar 2=1$. We can then recompute the connection as 
\begin{align*}\nabla s^i&=\theta\tens s^i-\sigma(s^i\tens\theta)={\imath \sigma_j\over 2}(s^j\tens s^i-\sigma(s^i\tens s^j)\\
&= {\imath\over 2}\sigma^i(s^1\tens s^1+s^2\tens s^2)+ {\imath\over 2}\eps_{ij}\sigma^j(s^2\tens s^1-s^1\tens s^2)-\rho\sigma^i s^{\bar i}\tens s^i,\end{align*}
which one can check agrees with what we had in terms of $s,t$. These coefficients define the $\Gamma^i{}_{jk}$. The inner product is $(s^2,s^2)=1=-(s^1,s^1)$ and zero off diagonal. 

The Riemann curvature and Ricci tensor for the  canonical lift $i(s\wedge t)={1\over 2}(s\tens t+t\tens s)$  translate to
\begin{align}\label{curvM2}
R_{\nabla}s^i=   -\imath(1+\rho^2 ){\rm Vol}\tens s^i,\quad {\rm Ricci}=-{1\over 2}(1+\rho^2)(s^1\tens s^1+s^2\tens s^2),
\end{align}
which is real in the sense of a metric for the $*$-structures, but not quantum symmetric. The Ricci scalar vanishes, so the QRG is Ricci-scalar flat but not necessarily flat. 

\subsection{Clifford action and spinor connection for alternate QRG} 

We next solve for $C^i$. Thus, the induced `Clifford relations' (\ref{cliffrel}) in the simplest case with $\varphi=\id$ and $\kappa=1$ come down in matrix terms to
\begin{equation}\label{Cirel} C^1 C^2= C^2 C^1=0,\quad -(C^1)^2+ (C^2)^2= 2{\id}\end{equation}
with the largest component (2-parameter family) of solutions 
\begin{equation}\label{cliffab} C^1=\pm\imath\begin{pmatrix}   \sqrt{2}- a &  b \\  \frac{a}{b}(\sqrt{2}- a) & a\end{pmatrix},\quad 
C^2=\pm\begin{pmatrix} a &  - b\\  -\frac{a}{b}(\sqrt{2}- a) &\sqrt{2}- a\end{pmatrix},\quad  a,b\in 
\C.\end{equation}
One can check that these $C^i$ admit no solutions for $J,\sigma_\CS$ and we must look more generally. The proof involves solving the covariance of $\la$ for parameters $a,b$, leading to a particular form of $\sigma_\CS$, which then fails (\ref{CJ}) using the fact that ${\rm Tr}(C^1)=\pm\imath\sqrt{2}$, ${\rm Tr}(C^2)=\pm\sqrt{2}$ are imaginary/real, which prescribes also the conjugated ${\rm Tr}(\bar C^i)$. We therefore look more generally by dropping the first half of the relations (\ref{Cirel}). 

\begin{proposition} If we impose only $-(C^1)^2+ (C^2)^2= 2\id$ then there are principally 5D and 4D moduli and we parametrise these by $c_{11},c_{12},d_{11},d_{12}$ and where applicable $d_{21}$ as three types of solutions,
\begin{enumerate}
\item[I:] $C^1=\left(
\begin{array}{cc}
 c_{11} & c_{12} \\
 \frac{-c_{11}^2+d_{11}^2+d_{12} d_{21}-2}{c_{12}} & -c_{11} \\
\end{array}
\right),\quad C^2=\left(
\begin{array}{cc}
 d_{11} & d_{12} \\
 d_{21} & -d_{11} \\
\end{array}
\right);$

\item[II:] $C^1=\left(
\begin{array}{cc}
 c_{11} & c_{12} \\
 \frac{c_{12} \left(c_{11}^2-d_{11}^2+2\right)}{d_{12}^2-c_{12}^2} & -c_{11} \\
\end{array}
\right),\quad C^2=\left(
\begin{array}{cc}
 d_{11} & d_{12} \\
 \frac{d_{12} \left(c_{11}^2-d_{11}^2+2\right)}{d_{12}^2-c_{12}^2} & -d_{11} \\
\end{array}
\right);$

\item[III:] $C^1=\left(
\begin{array}{cc}
 c_{11} & c_{12} \\
 \frac{c_{12} \left(c_{11}^2-d_{11}^2+2\right)}{d_{12}^2-c_{12}^2} & \frac{-2 c_{12} d_{11} d_{12}+c_{11} \left(c_{12}^2+d_{12}^2\right)}{d_{12}^2-c_{12}^2} \\
\end{array}
\right),$

\noindent$ C^2=\left(
\begin{array}{cc}
 d_{11} & d_{12} \\
 \frac{d_{12} \left(c_{11}^2-d_{11}^2+2\right)}{d_{12}^2-c_{12}^2} & \frac{-2 c_{11} c_{12} d_{12}+d_{11}(c_{12}^2 + d_{12}^2)}{c_{12}^2-d_{12}^2} \\
\end{array}
\right).$
\end{enumerate}
There are then no solutions  for $J,\sigma_\CS$ in type I or for type III in the case (\ref{cliffab}). 
\end{proposition}
\proof This is by calculations using Mathematica. In all types, solving for the covariance of $\la$ in Lemma~\ref{leminn} gives a 2-parameter moduli
of the $\sigma_\CS, \rho$ for generic $C^i$ in the class. For the last part, we note that (\ref{CJ}) says that $\bar C^i$ is conjugate to $C^i$, from which it follows that
\begin{equation}\label{traceid} {\rm Tr}(\sigma_\CS{}^i{}_jC^j)=0,\quad -(\sigma_\CS{}^1{}_iC^i)^2+ (\sigma_\CS{}^2{}_iC^i)^2=2\id\end{equation}
in type I,II for the trace identities and in all types for the other identity. In type I, we solve the covariance of $\la$ condition and the trace identities together. There are no solutions for general $C^i$ but if we solve also for one of their entries then we obtain 1-parameter solutions where $(\sigma_\CS{}^1{}_1)_{11}$ is undetermined but $d_{21}$ is prescribed. We are then not able to impose the other half of (\ref{traceid}). The analysis for general type III appears rather more complex as we have only the 2nd half of (\ref{traceid}), but (\ref{cliffab}) is a special case where the trace has a fixed property under complex conjugation as noted and similar methods then apply to establish no solutions. \endproof

The case of type II therefore appears most natural and we are led to focus on this. In this case, we are forced to 
\[\rho= -{c_{12}^2 + d_{12}^2\over 2 c_{12} d_{12}},\quad d_{11}= {c_{11}d_{12}\over c_{12}},\quad  C^1=\left(
\begin{array}{cc}
 c_{11} & c_{12} \\
 -\frac{c_{12} \left(c_{11}^2-d_{11}^2+2\right)}{c_{12}^2-d_{12}^2} & -c_{11} \\
\end{array}
\right),\quad C^2=\left(
\begin{array}{cc}
 d_{11} & d_{12} \\
 \frac{d_{12} \left(c_{11}^2-d_{11}^2+2\right)}{d_{12}^2-c_{12}^2} & -d_{11} \\
\end{array}
\right)\]
and a reduced 5-parameter moduli of the $\sigma_\CS$ entries as the solution of the covariance of $\la$ condition and both parts of (\ref{traceid}). We also need 
\begin{equation}\label{imrho} \bar\rho+\rho=0,\quad \overline{\det C^i}=\eps'\det(\sigma_\CS{}^i{}_jC^j)\end{equation}
as further constraints independently of $J$. 

\begin{lemma}\label{lemCi} For the type II Clifford algebra, the covariance of $\la$ condition, the condition (\ref{CJ}) and that $\rho$ is imaginary lead to one fewer parameter in $\sigma_\CS$ (so the latter drops to 4-dimensional), $\eps'=1$ and uniquely determined $\overline{C^i}$ in terms of $C^i$ for given $\eps$ for both types of $J$. \end{lemma}
\proof Building on the analysis from the previous lemma, we impose (\ref{imrho}) which fixes one of the coefficients of $\sigma_\CS$ and one of the 3 entries of $\bar C^i$ as 
\[ { \overline{d_{12}}\over d_{12}}=-{\overline{c_{12}}\over c_{12}}\] 
We then solve for  (\ref{CJ}) to find the remaining entries of $\bar C^i$. For both types of $J$, we need $\eps'=1$ and for type (1)  $J$, we find
\begin{align*}\overline{c_{12}}&=\frac{(c_{12}^2-d_{12}^2)(c_{12}z-c_{11})-2 c_{12}^2}{c_{12} \eps  \left(c_{12}^2-d_{12}^2\right)},\\
\overline{c_{11}}&=\frac{(c_{12}^2-d_{11}^2)(c_{12} z-c_{11} )( c_{12} ( \eps -z\bar z)+c_{11}  \bar z) -2 c_{12}^2 \bar z }{c_{12} \eps  \left(c_{12}^2-d_{12}^2\right)},\end{align*}
while for type (2) $J$,  the case of $\eps=1$ is simpler with
\[ -{\overline{c_{12}}\over c_{12}}= \frac{\bar z}{z},\quad\overline{c_{11}}-c_{11}=c_{12}\bar z. \]
These conditions are  solved by writing $z:=|z|e^{\imath \theta}$, then
\begin{equation}\label{J2realCi}  d_{12}=t e^{\imath \theta}={t z\over |z|},\quad c_{12}=\imath s e^{\imath\theta}={\imath s z\over |z|},\quad c_{11}=x-{\imath\over 2}|z| s,\quad  x,s,t\in \R.\end{equation}
The full $C^i$ matrices are then determined according to the above as
\begin{align}\label{solC1}C^1=\left(
\begin{array}{cc}
 x -\frac{i s  |z| }{2} &{ i s z\over |z|} \\
 -\frac{i \bar z \left((s^2+t ^2) (s|z|+2i x)^2-8s^2\right)}{4 s |z| \left(s ^2+t ^2\right)} & \frac{i s  |z| }{2}-x\\
\end{array}
\right),\\
\label{solC2} C^2=\left(
\begin{array}{cc}
 -\frac{t  (s  |z| +2 i x )}{2 s } & {t  z\over |z| } \\
 -\frac{t  \bar z\left((s^2+t ^2)(s|z|+2i x)^2-8s^2\right)}{4 |z| s ^2 \left(s ^2+t ^2\right)} & \frac{1}{2} t  \left(|z| +\frac{2 i x }{s }\right) \\
\end{array}
\right),\end{align}
and we also had to fix 
\begin{equation}\label{solrho} \rho={\imath\over 2}(\frac{t}{s}-\frac{s}{t}).\end{equation}
\endproof

This lemma gives us relations between the $c_{ij},d_{ij}$ and their conjugates, i.e. reality conditions in the $C^i$. It remains to solve the $\bar\sigma_\CS$ identity in Lemma~\ref{leminn} for some $J$. 
The second identity in Lemma~\ref{leminn}  can be seen as defining the conjugates of the entries of $\sigma_\CS$, 
\begin{equation}\label{CSJ} \overline{\sigma_\CS{}^i{}_j}= J\sigma_\CS^{-1}{}^i{}_jJ^{-1},\end{equation}
where the inverse is as an $M_2(\C)$-valued matrix with respect to the $i,j$ indices. This gives us relations between the entries of $\sigma_\CS$ and their conjugates.  We focus on the simplest case $J$ of type (2) with $\eps=1$ identified in the proof of the last lemma.  

\begin{lemma} For $J$ of type (2) with $\eps=1$ and the 3 real parameters $s,t,x$ of $C^i$ identified in (\ref{J2realCi}), (\ref{CSJ}) imposes two further relations of $\sigma_\CS$, so this drops down to two parameters subject to a reality constraint. Solutions of this include a natural  $\sigma_\CS$ determined uniquely by the parameters of $C^i$ and the complex parameter $z$ of $J$.
\end{lemma}
\proof In principle, the relations for $\overline{\sigma_\CS}$ entail in particular 
\[ \overline{{\rm Tr}\sigma_\CS{}^i{}_j}= {\rm Tr}( \sigma_\CS^{-1}{}^i{}_j),\quad \overline{\det\sigma_\CS{}^i{}_j} = (\det \sigma_\CS{})^{-1}{}^i{}_j,\]
where we trace and det over the suppressed spinor matrix indices. These are independent of $J$, hence it may be helpful to solve for  these first.  In practice, for our example, it is cleaner to look at all the 12 equations for the constrained entries of $\sigma_\CS$ (since there are four parameters left), take their complex conjugates and compare with $\overline{\sigma_\CS}$ as defined by the right hand side of (\ref{CSJ}), solving the simplest (and $J$-independent) ones first. The end result is to reduce our moduli space of $\sigma_\CS$ entries to two dimensions (we add two more relations) and find that then (\ref{CSJ}) is solved. We choose these parameters as
\[ \eta:=(\sigma_\CS{}^2{}_1){}^1{}_1,\quad \zeta:=(\sigma_\CS{}^2{}_1){}^1{}_2\]
and their reality condition is
\begin{align*}  \bar\eta&=\frac{d_{12} \left(4 \zeta ^2 c_{12}^2+ \left(c_{12}^2-d_{12}^2\right)(d_{12}(\eta  c_{12} -\zeta  2 c_{11}-\zeta  c_{12}  \bar z)+2(-\zeta c_{11}+\eta c_{12})^2 )\right)}{c_{12} \left(8 \zeta ^2  c_{12}^2+(c_{12}^2-d_{12}^2)(-2 c_{11} \zeta   +2 \eta  c_{12}+d_{12})^2 \right)}, \\
 \bar\zeta&=\frac{\zeta  d_{12}^2 \left(c_{12}^2-d_{12}^2\right)\bar z}{z \left(8 \zeta^2  c_{12}^2+(c_{12}^2-d_{12}^2)(-2 c_{11} \zeta   +2 \eta  c_{12}+d_{12})^2 \right)}.\end{align*}
Given the solution (\ref{J2realCi}) for the $C^i$, a natural solution of these reality conditions is
\[ \eta=-\frac{(s ^2+t^2)tx |z|}{2 \left((s ^2+t^2) x ^2+2 s ^2\right)},\quad \zeta=\frac{t  x  \left(s ^2+t ^2\right)}{(s ^2+t^2) x ^2+2 s ^2}e^{\imath \theta}.\]
This was obtained by making a polar expansion $\eta=ue^{\imath\phi}, \zeta=v e^{\imath\psi}$ then asking for the imaginary part of  $\bar\zeta/(ve^{-\imath\psi})$ to vanish (where $\bar\zeta$ is given by the right hand side of the expression above). This gives a simple equation for $v$ with two solutions if $s,t$ are generic. Taking the one stated (the other does not lead to a solution of the whole system), we ask for the imaginary part of $\bar\eta/(ue^{-\imath\phi})$ to vanish with $\bar\eta$ from the right hand side of the above. This gives a complicated equation for $u$ which, however, greatly simplifies if we assume $\psi=\theta$ and $\phi=0$. 

The entire $\sigma_\CS$ for this simplest choice of $\eta,\zeta$ is 
\begin{align}\label{sols11}  \sigma_\CS{}^1{}_1&=\left(
\begin{array}{cc}
 1-\frac{i s  x  |z|  \left(s ^2+t^2\right)}{2 \left((s^2+t^2) x ^2+2s^2\right)} & \frac{i z s  x  \left(s ^2+t^2\right)}{|z| \left((s^2+t  ^2) x^2+2s^2\right)} \\
- \frac{ i \bar z |z| s x  (s^2+t  ^2)  }{4 \left((s^2+t^2) x^2+2s^2\right)} +{\bar z\over 2}& \frac{i s  x  |z|  \left(s ^2+t  ^2\right)}{2 \left((s^2+t  ^2) x ^2+2s ^2\right)} \\
\end{array}
\right),\\
\label{sols12}  \sigma_\CS{}^1{}_2&=\left(
\begin{array}{cc}
 -\frac{s ^2 x  |z|  \left(s ^2+t  ^2\right)}{2t \left((s^2+t  ^2) x^2+2s^2\right)} & \frac{z s ^2  x  \left(s ^2+t  ^2\right)}{|z|t\left((s^2+t  ^2) x^2+2s^2\right)} \\
 -\frac{ \bar z  |z| s^2 x (s^2+t^2)}{4 t\left((s^2+t^2) x^2+2s^2\right)} -{i \bar z s\over 2t}& \frac{ s^2  x  |z|  \left(s ^2+t  ^2\right)}{2 t \left((s^2+t  ^2) x ^2+2s ^2\right)} +\frac{i s}{t}\\
\end{array}
\right), \\
\label{sols21}  \sigma_\CS{}^2{}_1&= \left(
\begin{array}{cc}
 -\frac{t  x  |z|  \left(s ^2+t ^2\right)}{2 \left((s^2+t  ^2) x^2+2s^2\right)} & \frac{z t x  \left(s ^2+t ^2\right)}{|z| (s^2+t^2) x^2+2s^2)} \\
 -\frac{t  \bar z |z|   x   (s^2+ t^2)}{4 \left((s^2+t  ^2) x^2+2s^2\right)} - {i \bar z t\over 2s}&  \frac{t  x   |z|  \left(s ^2+t ^2\right)}{2 (\left(s^2+t^2\right)x^2+2 s^2)}+\frac{i t }{s }\\
\end{array}
\right),\\
\label{sols22}   \sigma_\CS{}^2{}_2&=   \left(
\begin{array}{cc}
-1+\frac{i s x |z| \left(s^2+t^2\right)}{2 (\left(s^2+t^2\right)x^2 +2 s^2)}  & -\frac{i z s x  \left(s ^2+t  ^2\right)}{|z|\left((s^2+t  ^2) x^2+2s^2\right)} \\
 \frac{i \bar z |z|  s x  \left(s^2+t^2\right)}{4 \left((s^2+t^2) x^2+2s^2\right)} -{ \bar z\over 2} & -\frac{i s  x  |z|  \left(s^2+t^2\right)}{2 \left((s^2+t^2) x ^2+2s^2\right)}\\
\end{array}
\right). \end{align}
 \endproof
 
\subsection{Construction of non-hermitian almost spectral triples on $M_2(\C)$}

We next consider the Hilbert space structure, taking the obvious choice
\begin{equation}\label{hilbM2} \<\phi_\alpha e^\alpha,\psi_\beta e^\beta\>={1\over 2}{\rm Tr}(\phi_\alpha^*\psi_\alpha),\end{equation}
corresponding  to $\int={1\over 2}{\rm Tr}$ and $\mu^{\alpha\beta}=\delta^{\alpha\beta}$ in (\ref{hilbS}). Here, the spinor coefficients $\phi_\alpha,\psi_\alpha\in M_2(\C)$ with its usual (hermitian adjoint) $*$-operation.
 
 \begin{theorem}\label{M2evenspec} On $M_2(\C)$ with its alternate 1-parameter quantum Riemannian geometry as in Section~\ref{secM2qrgalt}, the conditions for a geometrically realised spectral triple at the local tensorial level have a natural solution with 4 real parameters $s,t,x,y$, $\eps=\eps'=\eps''=1$,
 \[ J=\begin{pmatrix} 1 & 0\\ 0 & e^{2\imath y}\end{pmatrix},\quad \gamma=\pm\begin{pmatrix} 1 & 2 { e^{\imath y} s x(s^2+t^2)\over (s^2+t^2)x^2+2 s^2}\\ 0 & -1\end{pmatrix}\]
 such that $\CJ$ is an antilinear isometry with respect to (\ref{hilbM2}). Moreover, the QLC has its parameter $\rho$ fixed as (\ref{solrho}) with curvature proportional to
 \[1+\rho^2=\frac{1}{4} \left(6-\frac{s^2}{t^2}-\frac{t^2}{s^2}\right).\]
 \end{theorem} 
 \proof We first bring together our lemmas above. We fix $J$ of type (2) with $\eps=1$ as this gave the most natural solutions for $C^i$. Here $z$ is a free parameter in $J$ in case (2) of Lemma~\ref{lemJ}. The three parameters $s,t,x$ go into solving for $C^i$ to obey the 2nd half of (\ref{Cirel}), covariance of $\la$ and (\ref{CJ}), and also determine $\rho$, i.e. the choice of QLC which was otherwise not fixed for the given metric. Existence of the $C^i$ also cuts down the allowed $\nabla_\CS$ to 4 parameters and among solutions we found a natural choice of these so as to obey the reality condition in Lemma~\ref{leminn} for $\nabla_\CS$ to be $*$-preserving. 
 
 Next, we take the natural choice of Hilbert space structure on $\CS=M_2(\C)\tens\C^2$ as stated. As $\eps=1$, we then need $J$ to be symmetric. We achieve this by taking a limiting value $z\to 0$ along a ray such that $-z/\bar z= e^{2\imath y}$ for some real $y$. Our solutions for $C^i,\sigma_\CS{}^i{}_j$ then simplify considerably to
 \begin{equation}\label{finalC}C^1=\left(
\begin{array}{cc}
 x &s e^{\imath y} \\
- \frac{(s^2+t^2)x^2+2s^2}{ e^{\imath y} s  \left(s ^2+t ^2\right)} & -x\\
\end{array}
\right),\quad  C^2=-{\imath t\over s} C^1,\end{equation}
\begin{equation}\label{finals}  \sigma_\CS{}^1{}_1=\left(
\begin{array}{cc}
 1 & \frac{e^{\imath y} s  x  \left(s ^2+t^2\right)}{  (s^2+t^2)x^2+2s^2} \\
0& 0 \\
\end{array}
\right) =-\sigma_\CS{}^2{}_2,\quad  \sigma_\CS{}^1{}_2={\imath s\over t}\left(
\begin{array}{cc}
0 &-\frac{ e^{\imath y}s   x  \left(s ^2+t  ^2\right)}{ (s^2+t^2)x^2 +2s^2} \\
0& 1\\
\end{array}
\right) ={s^2\over t^2}\sigma_\CS{}^2{}_1.  \end{equation}

For the last part, we first solve (\ref{gamsq}),(\ref{Cgam}) and (\ref{Sgaminn}) to find the unique answer for $\gamma$ up to sign (which is not determined by the equations). Then (\ref{Jgam}) is solved automatically with $\eps''=1$. \endproof

Note that we can set $y=0$ by conjugation by a unitary ${\rm diag}(1,e^{\imath y})$. We also see that there are four points $t=\pm(\sqrt{2}\pm' 1)s$ where the QRG is flat. It remains to see if $\dirac$ and $\gamma$ constructed above are (anti)hermitian. 

\begin{proposition}\label{propM2dirac} The geometrically realised Dirac operator in Theorem~\ref{M2evenspec} has the form
\[ (\dirac \psi)_\gamma={\imath\over 2}([\sigma^1,\psi_\alpha] C^{1\alpha}{}_\gamma+\{\sigma^2,\psi_\alpha\} C^{2\alpha}{}_\gamma )\]
which we write as
\[ \dirac=\slashed{\del}^{(1)}- {\imath t\over s}\slashed{\del}^{(2)},\quad (\slashed{\del}^{(1)}\psi)_\gamma={\imath\over 2}[\sigma^1,\psi_\alpha]C^{1\alpha}{}_\gamma,\quad (\slashed{\del}^{(2)}\psi)_\gamma={\imath\over 2}\{\sigma^2,\psi_\alpha\} C^{1\alpha}{}_\gamma.\]
If  $x=\pm s\sqrt{1-{2\over s^2+t^2}}$ then $\slashed{\del}^{(i)}$ are both hermitian with respect to (\ref{hilbM2}), while $\gamma$ is hermitian if $x=0$. All three are hermitian if $s^2+t^2=2$ and $x=0$.  \end{proposition} 
 \proof  We again use (\ref{diracinner}) given the inner form of the connection and the form of $\theta$ for this calculus. This time, however,  for  our particular 4-parameter solutions, we have
 \[ \sigma_\CS{}^j{}_i C^i=\pm C^j\]
  with +1 if $j=1$ and -1 if $j=2$, which gives the result stated. Next, if $C^1$ is {\em antihermitian}, which happens for $x$ with the value stated, then 
  \[\<\slashed{\del}^{(1)}\psi,\phi\>-\<\psi,\slashed{\del}^{(1)}\phi\>={\imath C^{1\gamma}{}_\alpha \over 4}{\rm Tr}\left( [
  \sigma^1,\psi_\alpha]^\dagger \phi_\gamma- \psi_\alpha^\dagger[\sigma^1,\phi_\gamma]\right)=0,\]
   \[\<\slashed{\del}^{(2)}\psi,\phi\>-\<\psi,\slashed{\del}^{(2)}\phi\>={\imath C^{1\gamma}{}_\alpha \over 4}{\rm Tr}\left( \{
  \sigma^2,\psi_\alpha\}^\dagger \phi_\gamma- \psi_\alpha^\dagger\{\sigma^2,\phi_\gamma\}\right)=0,\]
  using the trace property.

 For $\gamma$ to be hermitian with respect to (\ref{hilbM2}) amount to the matrix for $\gamma$ hermitian, which is clearly not the case unless $x s=0$. Our solutions suppose $s\ne 0$ so the condition is $x=0$ and then $\gamma=\pm\sigma^3$. This is compatible with our condition for $\slashed{\del}^{(i)}$ when $s^2+t^2=2$. 
  \endproof
 
We see that  our 4-parameter space of solutions in Theorem~\ref{M2evenspec}  reduces to a circle $s^2+t^2=2$ (with $s,t\ne 0$) if we require $\slashed{\del}^{(i)},\gamma$ hermitian. Choosing the $+$ sign for square roots and setting $y=0$ by the freedom under unitary equivalence, we have in this case
 \begin{align}   C^1&=\imath s\sigma^2   ,\quad  C^2=  t \sigma^2,\quad  J=\id,\\
 \sigma_\CS{}^1{}_1&=E_{11}=-\sigma_\CS{}^2{}_2 ,\quad \sigma_\CS{}^1{}_2={\imath s\over t}E_{22},\quad \sigma_\CS{}^2{}_1={\imath t\over s}E_{22} ,\quad \gamma=\sigma^3.\end{align}
  In fact, the pattern of signs would fit into the period 8 KO-dimension signs with $n=0$ if $\dirac$ were to be hermitian and would fit $n=1$ if $\dirac$ were to be antihermitian (then we would flip the value of $\eps'$), but our $\dirac$ is neither: the first part of $\dirac$ is hermitian and the second part is antihermitian. Similarly, $C^1$ is antihermitian and $C^2$ is hermitian in line with the Lorentzian nature of the metric  (\ref{M2g_alt}) on $M_2(\C)$.

 We do not claim that Theorem~\ref{M2evenspec}  is the only data at the local tensorial level for a geometrically realised spectral triples on $M_2(\C)$ for this quantum metric, but our results were obtained by solving the equations with the simplest nontrivial choices at the various branches as per the derivation above.   Also note that by Lemma~\ref{leminn}, we are free (both here and for the metric in the main text) to add a bimodule map in the form of matrices $A_i$ to the connection $\nabla_\CS$ subject to the conditions in the lemma.  This will appear as further freedom in the construction of a geometrically realised Dirac operator as fluctuations of the inner case that we solved for. On the other hand, it does not appear that looking more widely will solve the problem that $\dirac$ is not antihermitian for the quantum metric in this appendix.
  
\section*{Declaration}

Data sharing is not applicable as no datasets were generated or analysed during the current study. The authors have no competing interests to declare that are relevant to the content of this article.


\begin{thebibliography}{ggghhh}

\bibitem{ArgMa4}J. Argota-Quiroz and S. Majid, Quantum gravity on finite spacetimes and dynamical mass, 
PoS (2022) 210 (41pp)

 \bibitem{BegMa} E.J. Beggs and S. Majid, {\em Quantum Riemannian Geometry},  Grundlehren der mathematischen Wissenschaften, Vol. 355, Springer (2020) 809pp
  
\bibitem{BG}J.W. Barrett and J. Gaunt,  Finite spectral triples for the fuzzy torus, arXiv:1908.06796 (math.QA)
 
 \bibitem{BegMa:twi} E.J. Beggs and S. Majid, Nonassociative Riemannian geometry by twisting,  J. Phys. Conf. Ser. 254 (2010) 012002 (29pp)
 
 \bibitem{BegMa:spe}E.J. Beggs and S. Majid, Spectral triples from bimodule connections and Chern connections, J. Noncomm. Geom., 11 (2017) 669--701
 
 \bibitem{BegMa:cur}E. Beggs and S. Majid, Quantum geodesics and curvature, in press Lett. Math. Phys. (2023)
 
 \bibitem{CPR}A.L. Carey, J. Phillips and A. Rennie, Spectral triples: examples and index theory, in {\em Noncommutative Geometry and Physics: Renormalisation, Motives, Index Theory}  (European Mathematical Society: Zurich) 2011, pp175--265
 
 \bibitem{Con:dif}A. Connes, C*-alg\`ebres et g\'eom\'etrie differentielle, C.R. Acad. Sci. Paris, Ser. A-B , 290, (1980) 
\bibitem{Con0} A. Connes, Gravity coupled with matter and foundation of noncommutative geometry, Commun. Math. Phys. 182 (1996) 155

\bibitem{Con}
A. Connes,   Noncommutative Geometry, 
Academic Press, Inc., San Diego, CA, 1994

\bibitem{ConMar} A. Connes and M. Marcolli, {\em Noncommutative Geometry, Quantum Fields and Motives} (AMS Colloquium Publications Vol 55), Hindustan Book Agency, 2008. 

\bibitem{DabSit}L. Dabrowski and A. Sitarz, Fermion masses, mass-mixing and the almost commutative geometry of the Standard Model,  JHEP (2019) 68 

\bibitem{And}F. D'Andrea, F.  Lizzi and J. C. Varilly, Metric properties of the fuzzy sphere, Lett. Math. Phys. 103 (2013) 183--205

\bibitem{Dev} A. Devastato, S. Farnsworth, F. Lizzi and P. Martinetti, Lorentz signature and twisted spectral triples, JHEP  (2018) 89


\bibitem{DFR}S. Doplicher, K. Fredenhagen and J. E. Roberts, The quantum structure of spacetime at the Planck scale and quantum fields, Commun. Math. Phys. 172 (1995) 187--220


\bibitem{DVMic}
M. Dubois-Violette and  P.W.\ Michor, Connections on central bimodules in 
noncommutative differential geometry, J.\ Geom.\ Phys.\ 20 (1996) 218--232

\bibitem{Bon}J. M. Gracia-Bond\'ia, J.C. V\'arilly, H. Figueroa,  {\em Elements of Noncommutative Geometry}, Birkh\"auser,  Boston, 2001

\bibitem{Hoo}G. 't Hooft, Quantization of point particles in 2+1 dimensional gravity and space- time discreteness, Class. Quant.  Grav. 13 (1996) 1023

\bibitem{LirMa}E. Lira-Torres and S. Majid, Quantum gravity and Riemannian geometry on the fuzzy sphere, Lett. Math. Phys. (2021) 111:29 (21pp)

\bibitem{LirMa2}E. Lira-Torres and S. Majid, Geometric Dirac operator on the fuzzy sphere,  Lett. Math. Phys. (2022) 112:10


\bibitem{LiuMa2} C. Liu and S. Majid, Quantum Kaluza-Klein theory with $M_2(\C)$, arXiv: 2303.06239 (gr-qc)



\bibitem{Ma:pla}S. Majid, Hopf algebras for physics at the Planck scale, Class. Quant.  Grav. 5 (1988) 1587--1607

\bibitem{MaRue}S. Majid and H. Ruegg, Bicrossproduct structure of the $\kappa$-Poincare group and non-commutative geometry, Phys. Lett. B. 334 (1994) 348--354


\bibitem{Ma:dir}S. Majid, Dirac operator associated to a quantum metric, arXiv: 2302.05891 (math.QA)

\bibitem{MaSim} S. Majid and F. Simao, Quantum jet bundles, arXiv: 2202.03067 (math.QA)

\bibitem{Mou}
J.\ Mourad, Linear connections in noncommutative geometry, Class.\ Quant.   
Grav.\ 12 (1995)  965--974


\bibitem{PasSit1}M. Paschke and A.Sitarz, On spin structures and Dirac operators on the noncommutative torus, Lett. Math. Phys. 77 (2006), 317--327

\bibitem{PasSit} M. Paschke and A. Sitarz, Equivariant Lorentzian spectral triples, arXiv:0611029 (math-ph)


\bibitem{Sny}H.S. Snyder,  Quantized space-time, Phys. Rev. 71 (1947) 38--41


\bibitem{Ven}J.J. Venselaar, Classification of spin structures on the noncommutative n-torus, J. Noncomm. Geom. 7 (2013), 787--816

 \end{thebibliography}
\end{document}